\documentclass[12pt]{amsart}

\usepackage{amssymb,geometry,stmaryrd,color}

\usepackage{graphicx}
\usepackage{caption}
\begin{document}
\baselineskip=18pt
\setcounter{page}{1}
    
\newtheorem{conjecture}{Conjecture}
\newtheorem{corollary}{Corollary}
\newtheorem{DEF}{Definition\!\!}
\newtheorem{theorem}{Theorem}
\newtheorem{lemma}{Lemma}
\newtheorem{proposition}{Proposition}
\newtheorem{remark}{Remark}

\renewcommand{\theDEF}{}

\def\a{\alpha}
\def\b{\beta}
\def\B{{\bf B}} 
\def\C{{\bf C}} 
\def\K{{\bf K}}
\def\BB{{\mathcal{B}}} 
\def\DD{{\mathcal{D}}} 
\def\GG{{\mathcal{G}}} 
\def\KK{{\mathcal{K}}} 
\def\LL{{\mathcal{L}}} 
\def\SS{{\mathcal{S}}}
\def\UU{{\mathcal{U}}}
\def\ca{c_{\a}}
\def\ka{\kappa_{\a}}
\def\coa{c_{\a, 0}}
\def\cua{c_{\a, u}}
\def\cL{{\mathcal{L}}} 
\def\cV{{\mathcal{V}}} 
\def\Ea{E_\a}
\def\eps{{\varepsilon}} 
\def\esp{{\mathbb{E}}} 
\def\Ga{{\Gamma}} 
\def\G{{\bf \Gamma}} 
\def\e{{\rm e}}
\def\i{{\rm i}}
\def\L{{\bf L}}
\def\lbd{\lambda}
\def\lacc{\left\{}
\def\lcr{\left[}
\def\lpa{\left(}
\def\lva{\left|}
\def\M{{\bf M}}
\def\Nn{{\bf N}}
\def\prst{{\leq_{st}}}
\def\prost{{\prec_{st}}}
\def\prcvx{{\prec_{cx}}}
\def\Rr{{\bf R}}
\def\CC{{\mathbb{C}}}
\def\NN{{\mathbb{N}}} 
\def\QQ{{\mathbb{Q}}} 
\def\pb{{\mathbb{P}}}
\def\rl{{\mathbb{R}}}
\def\rpa{\right)}
\def\rcr{\right]}
\def\rva{\right|}
\def\Tt{{\bf \Theta}}
\def\Ttt{{\tilde \Tt}}
\def\W{{\bf W}}
\def\X{{\bf X}}
\def\XX{{\mathcal X}}
\def\YY{{\mathcal Y}}
\def\U{{\bf U}}
\def\V{{\bf V}}
\def\Un{{\bf 1}}
\def\Y{{\bf Y}}
\def\Z{{\bf Z}}
\def\ZZ{{\mathbb{Z}}}
\def\A{{\bf A}}
\def\AA{{\mathcal A}}
\def\hAA{{\hat \AA}}
\def\hL{{\hat L}}
\def\hT{{\hat T}}

\def\claw{\stackrel{d}{\longrightarrow}}
\def\elaw{\stackrel{d}{=}}
\def\pslaw{\stackrel{a.s.}{\longrightarrow}}
\def\qed{\hfill$\square$}

\newcommand*\pFqskip{8mu}
\catcode`,\active
\newcommand*\pFq{\begingroup
        \catcode`\,\active
        \def ,{\mskip\pFqskip\relax}%
        \dopFq
}
\catcode`\,12
\def\dopFq#1#2#3#4#5{%
        {}_{#1}F_{#2}\biggl[\genfrac..{0pt}{}{#3}{#4};#5\biggr]%
        \endgroup
}

%%%%%%%% Takahiro's macros %%%%%%%%%

\newcommand{\R}{\mathbb{R}}
\newcommand{\N}{\mathbb{N}}
\newcommand{\bP}{\mathbb{P}}    % Probability 
\newcommand{\bE}{\mathbb{E}}    % Expectation
\def\ii{{\rm i}}
\def\S{\mathbf{S}}
\def\F{\mathbf{T}}
\def\W{\mathbf{W}}

\title[Some properties of the free stable distributions]{Some properties of the free stable distributions}

\author[Takahiro Hasebe]{Takahiro Hasebe}

\address{Department of Mathematics, Hokkaido University, North 10 West 8, Kita-Ku, Sapporo 060-0810, Japan. {\em Email}: {\tt thasebe@math.sci.hokudai.ac.jp}}

\author[Thomas Simon]{Thomas Simon}

\address{Laboratoire Paul Painlev\'e, Universit\'e de Lille, Cit\'e Scientifique, 59655 Villeneuve d'Ascq Cedex, France. {\em Email}: {\tt simon@math.univ-lille1.fr}}

\author[Min Wang]{Min Wang}

\address{Laboratoire Paul Painlev\'e, Universit\'e de Lille, Cit\'e Scientifique, 59655 Villeneuve d'Ascq Cedex, France. {\em Email}: {\tt min.wang@math.univ-lille1.fr}}



\begin{abstract} We investigate certain analytical properties of the free $\a-$stable densities on the line. We prove that they are all classically infinitely divisible when $\a\le 1$, and that they belong to the extended Thorin class when $\a\le 3/4.$ The L\'evy measure is explicitly computed for $\a =1,$ showing that the free 1-stable random variables are not Thorin except in the drifted Cauchy case. In the symmetric case we show that the free stable densities are not infinitely divisible when $\a > 1.$ In the one-sided case we prove, refining unimodality, that the densities are whale-shaped that is their successive derivatives vanish exactly once. Finally, we derive a collection of results connected to the fine structure of the one-sided free stable densities, including a detailed analysis of the Kanter random variable, complete asymptotic expansions at zero, a new identity for the Beta-Gamma algebra, and several intrinsic properties of whale-shaped densities.
\end{abstract}

\maketitle

\section{Introduction}

In this paper, we investigate certain properties of free stable densities on the line. The latter are the solutions $f$ to the following convolution equation
\begin{equation}
\label{Main}
a\X_1\; +\; b\X_2 \; \elaw \; c \X\; +\; d
\end{equation}
where $\X_1, \X_2$ are free independent copies of a random variable $\X$ with density $f$, $a,b$ are arbitrary positive real numbers, $c$ is a positive real number depending on $a, b$, and $d$ is a real number. As in the classical framework, it turns out that there exist solutions to (\ref{Main}) only if $c = (a^\a+b^\a)^{1/\a}$ for some fixed $\a\in (0,2].$ We will be mostly concerned with free strictly stable densities, which correspond to the case $d=0.$ In this framework the Voiculescu transform of $f$ writes, up to multiplicative normalization,
\begin{equation}
\label{Voicu}
\phi_{\a,\rho} (z) \; =\; - e^{\i \pi \a\rho} z^{-\a +1}, \quad \Im (z) > 0,
\end{equation}
with $\rho \in [0,1]$ if $\a\in (0,1]$ and $\rho\in[1-1/\a, 1/\a]$ if $\a\in[1,2].$ We refer e.g. to \cite{NS} for some background on the free additive convolution, to \cite{BerVoi} for the original solution to the equation (\ref{Main}), and to the introduction of \cite{HK} for the above parametrization $(\a,\rho)$, which mimics that of the strict classical framework. Let us also recall that free stable laws appear as limit distributions of spectra of large random matrices with possibly unbounded variance - see \cite{BG,CD}, and that their domains of attraction have been fully characterized in \cite{BP95, BerPat}. In the following, we will denote by $\X_{\a,\rho}$ the random variable whose Voiculescu transform is given by (\ref{Voicu}), and set $f_{\a,\rho}$ for its density. The analogy with the classical case extends to the fact, observed in Corollary 1.3 of \cite{HK}, that with our parametrization one has
$$\pb[\X_{\a,\rho}\ge 0]\; =\; \rho.$$
 The explicit form of the Voiculescu transform also shows that $\X_{\a,\rho}\elaw - \X_{\a, 1-\rho}.$ In this paper, some focus will put on the one-sided case and we will use the shorter notations $\X_{\a,1} = \X_\a$ and $f_{\a,1} = f_\a.$ Throughout, the random variable $\X_{\a,\rho}$ will be mostly handled as a classical random variable via its usual Fourier, Laplace and Mellin transforms, except for a few situations where the free independence is discussed.

Several analytical properties of free stable densities have been derived in the Appendix to \cite{BerPat}, where it was shown in particular that they can be expressed in closed form via the inverse of certain trigonometric functions. It is also a consequence of Proposition 5.12 in \cite{BerPat} that save for $\a =1,$ every free $\a-$stable density is, as in the classical framework, an affine transformation of some $f_{\a,\rho}.$ The density $f_{\a,\rho}$ turns out to be a truly explicit function in three specific situations only, which is again reminiscent of the classical case:\\

\begin{itemize}

\item $f_{2, 1/2} (x) = {\displaystyle \frac{\sqrt{4-x^2}}{2\pi}} \;\;\mbox{for}\; x\in[-2,2],\;\; $ (semi-circular density),\\

\item $f_{1/2} (x) = {\displaystyle \frac{\sqrt{4x -1}}{2\pi x^2}} \;\;\mbox{for}\; x\ge 1/4,\;\;$ (inverse Beta density),\\

\item $f_{1, \rho} (x) = {\displaystyle \frac{\sin (\pi\rho)}{\pi (x^2 +2\cos (\pi\rho) x +1)}} \;\;\mbox{for}\; x\in \rl,\;\; $ (standard Cauchy density with drift).\\

\end{itemize}

The study of $f_{\a,\rho}$ was carried on further in \cite{HM, HK} where, among other results, several factorizations and series representations were obtained. Our purpose in this paper is to deduce from these results several new and non-trivial properties. Our first findings deal with the infinite divisibility of $\X_{\a,\rho}$. Since this random variable is freely infinitely divisible (FID), it is a natural question whether it is also classically infinitely divisible (ID). 

\begin{theorem}
\label{freeID} One has

{\em (a)} For every $\alpha \in (0,1]$ and $\rho \in [0,1]$, the random variable $\X_{\a, \rho}$ is {\em ID}.

{\em (b)} For every $\a \in (1,2],$ the random variable $\X_{\a, 1/2}$ is not {\em ID}.
\end{theorem}

Above, the non ID character of $\X_{2,1/2}$ is plain from the compactness of its support. Observe also that by continuity of the law of $\X_{\a,\rho}$ in $(\a,\rho)$ and closedness in law of the ID property - see e.g. Lemma 7.8 in \cite{Sato}, for every $\a\in(1,2)$ there exists some $\epsilon(\a) >0$ such that $\X_{\a,\rho}$ is not ID for all $\rho\in[1/2-\epsilon(\a),1/2 +\varepsilon(\a)].$ We believe that one can take $\varepsilon(\a) =1/\a -1/2,$ that is our above result is optimal with respect to the ID property. Unfortunately, we found no evidence for this fact as yet - see Remark \ref{IDtwo} for possible approaches. 

As it will turn out in the proof, for $\a\le 1$ the ID random variables $\X_{\a, \rho}$ have no Gaussian component. A natural question is then the structure of their L\'evy measure. We will say that the law of a positive ID random variable is a generalized Gamma convolution (GGC) if its L\'evy measure has a density $\varphi$ such that $x\varphi(x)$ is a completely monotonic (CM) function on $(0,+\infty).$ There exists an extensive literature on such positive distributions, starting from the seventies with the works of O. Thorin. The denomination comes from the fact that up to translation, these laws are those of the random integrals
$$\int_0^\infty a(t)\,d\G_t$$
where $a(t)$ is a suitable deterministic function and $\{\G_t, \, t\ge 0\}$ is the Gamma subordinator. We refer to \cite{B} for a comprehensive monograph with an accent on the Pick functions representation and to the more recent survey \cite{JRY} for the above Wiener-Gamma integral representation, among other topics. See also Chapters 8 and 9 in \cite{SSV} for their relationship with Stieltjes functions. In Chapter 7 of \cite{B}, this notion is extended to distributions on the real line. Following (7.1.5) therein, we will say that the law of a real ID random variable is an extended GGC if its L\'evy measure has a density $\varphi$ such that $x \varphi(x)$ and $x \varphi(-x)$ are CM as a function of $x$ on $(0,+\infty).$ In order to simplify our presentation, we will also use the notation GGC for extended GGC.
 
\begin{theorem}
\label{GGCID}
For every $\alpha\in (0,3/4]$ and $\rho\in [0,1]$, the law of $\X_{\a, \rho}$ is a {\em GGC.}
\end{theorem}

Contrary to the above, we think that this result is not optimal and that the random variable $\X_{\a, \rho}$ has a GGC law at least for every $\alpha\in (0,4/5]$ and $\rho\in [0,1]$ - see Conjecture \ref{KaGGC}. During our proof, we will show that  for every $\a,\rho < 1$ the GGC character of $\X_{\a,\rho}$ is a consequence of that of $\X_{\a}.$ Unfortunately this simpler question, which is connected to the hyperbolically completely monotonic (HCM) character of negative powers of the classical positive stable distribution, is rather involved. Moreover, we will see in Corollary \ref{NoGG} that the law of $\X_\a$ is not a GGC for $\a$ close enough to 1.\\

Our next result deals with the case $\a = 1.$ According to the Appendix of \cite{BerPat}, the Voiculescu transform writes here, up to affine transformation,
$$\phi_\rho (z) \; =\; - \frac{2}{\pi}\,(\rho\pi \i\; +\; (1-2\rho) \log z), \quad \Im (z) > 0,$$
for some $\rho\in [0,1].$ By (\ref{Voicu}), this means that a free 1-stable distribution is up to translation the law of the free independent sum 
$$\C_{a,b}\; \elaw\; a\X_{1,1/2}\;+\; b\F,$$ 
for $a\ge 0, b\in\rl,$ where $\F$ has Voiculescu transform $-\log z$ and will be called henceforth the exceptional free 1-stable random variable. For example, $\phi_{1/2}$ is the Voiculescu transform of $\X_{1,1/2},$ whereas $\phi_0$ is that of $\frac{2}{\pi}(\F + \log (\pi/2))$ and $\phi_1$ that of $-\frac{2}{\pi}(\F + \log (\pi/2)).$ The density of $\C_{a,b}$ can be retrieved from Proposition A.1.3 of \cite{BerPat}, in an implicit way. In this paper, taking advantage of a factorization due to Zolotarev for the exceptional classical 1-stable random variable, we obtain the following explicit result.

\begin{theorem}
\label{Free1}
The random variable $\C_{a,b}$ is {\em ID} without Gaussian component and with L\'evy measure $$\frac{1}{x^2}\lpa \frac{a}{\pi}\,\Un_{\{x\neq 0\}}\; +\; \vert b\vert \lpa 1 - \frac{\vert b^{-1}x\vert\, e^{-2\vert b^{-1}x \vert}}{1- e^{-\vert b^{-1}x\vert}}\rpa\Un_{\{ bx < 0\}}\rpa dx,$$
where the second term is assumed to be zero if $b=0$. 
\end{theorem}

This computation implies - see Remark \ref{MellWbis} - that the random variable $\C_{a,b}$ is self-decomposable (SD) and has CM jumps, but that its law is not a GGC except for $b = 0.$ A key-tool for the proof is an identity connecting $\F$ and the free Gumbel distribution - see Proposition \ref{prop:SF}, providing an analogue of Zolotarev's factorization in the free setting, and which is interesting in its own right.\\
        
Our last main result concerns the shape of the densities $f_{\a,\rho}.$ It was shown in the Appendix to \cite{BerPat} that the latter are analytic on the interior of their support, and strictly unimodal i.e. they have a unique local maximum. These basic properties mimic those of the classical stable densities displayed in the monograph \cite{Z}. A refinement of strict unimodality was recently investigated in \cite{Kw, TSPAMS}, where it is shown that the classical stable densities are bell-shaped (BS), that is their $n-$th derivative vanishes exactly $n$ times on the interior of their support, as is the case for the standard Gaussian density. The free strictly 1-stable density $f_{1,\rho}$ is BS, but it is visually clear that this property is not fulfilled neither by $f_{2,1/2}$ nor by $f_{1/2}.$ Let us introduce the following alternative refinement of strict unimodality.

\begin{DEF}
\label{whales}
A smooth non-negative function on $\rl$ is said to be {\em whale-shaped} if its support is a closed half-line, if it vanishes at both ends of its support, and if
$$\sharp\{ x\in {\rm Supp}\, f,\; f^{(n)}(x) = 0\}\; =\; 1$$
for every $n\ge 1.$
\end{DEF}

\begin{figure}
\begin{center}
%\advance\leftskip-10cm
%\advance\rightskip-7cm
\begin{minipage}{0.45\hsize}
\begin{center}
\includegraphics[scale=0.2, bb=0 0 960 731]{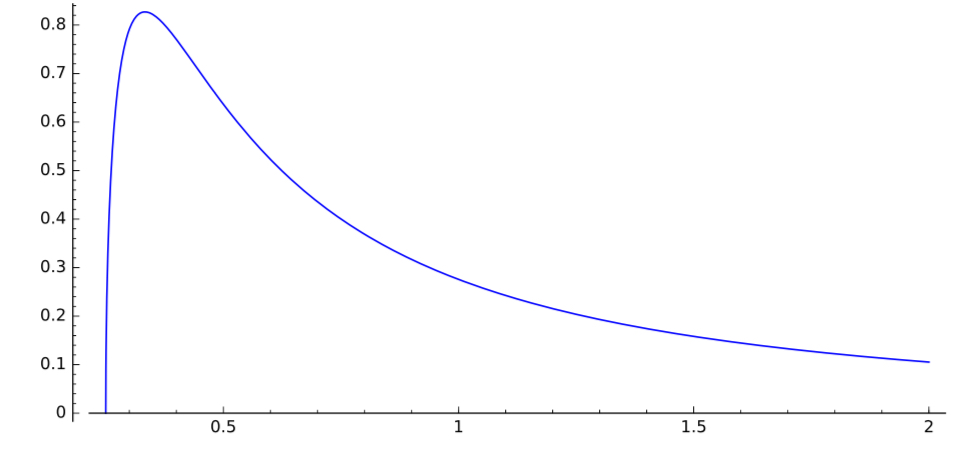} 
\caption{The free positive $1/2$-stable density (WS).} \label{WS1}
\end{center}
  \end{minipage}
\begin{minipage}{0.45\hsize}
\begin{center}
\includegraphics[scale=0.2, bb=0 0 960 731]{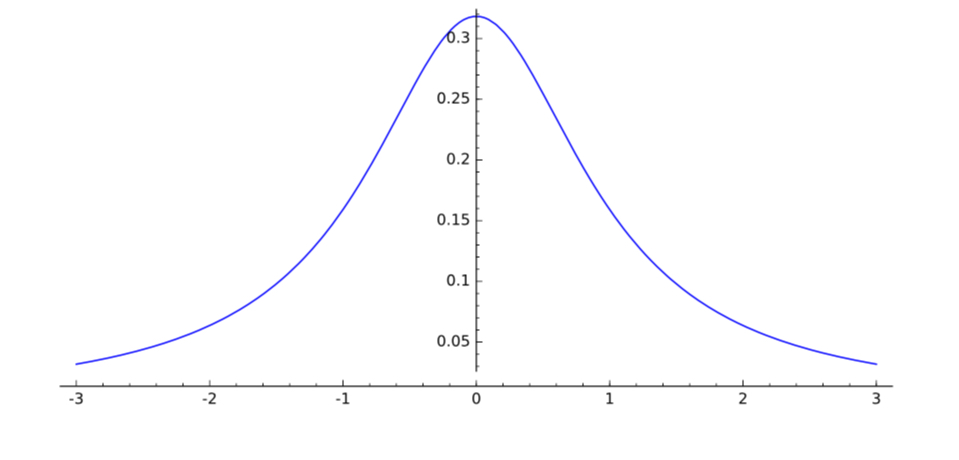}  
\caption{The free symmetric $1$-stable density (BS).} \label{BS1}
\end{center}
\end{minipage}
\end{center}
\end{figure}

The denomination comes from the visual aspect of such functions - see Figure \ref{WS1} and compare with the visual aspect of a bell-shaped density given in Figure \ref{BS1}. We will denote by WS the whale-shaped property and set ${\rm WS}_+$ (resp. ${\rm WS}_-$) for those whale-shaped functions whose support is a positive half-line $[x_0, +\infty)$ for some $x_0\in\rl,$ resp. a negative half-line $(-\infty, x_0].$ Observe that if $f\in{\rm WS}_+,$ then $x\mapsto f(-x)$ belongs to ${\rm WS}_-.$ It is easy to see that if $f\in{\rm WS}_+$ has support $[x_0, +\infty),$ then $f$ is positive on $(x_0,+\infty), f^{(n)}(+\infty) = 0$ and $(-1)^{n-1}f^{(n)}(x_0+) > 0$ for every $n\ge 1.$ In particular, the class ${\rm WBS}_0$ introduced in the main definition of \cite{TSPAMS} corresponds to those ${\rm WS}_+$ functions whose support is $(0,+\infty).$ Observe finally that the sequence of vanishing places of the successive derivatives of a function in ${\rm WS}_+$ increases, by Rolle's theorem. Other, less immediate, interesting properties of WS functions will be established in Section \ref{sec:WS}. 

\begin{theorem}
\label{shapes}
One has

{\em (a)} For every $\alpha\in (0,1),$ the density $f_\a$ is ${\rm WS}_+$

{\em (b)} For every $\a\in (0,3/4]$ and $\rho\in (0,1),$ the density $f_{\a,\rho}$ is {\em BS}

{\em (c)} The density of $\F$ is ${\rm WS}_-$

{\em (d)} For $a\neq 0$ and for $b=0$ or $ab^{-1}\in\pi\ZZ,$ the density of $\C_{a,b}$ is ${\rm BS}.$  

\end{theorem}

This result leaves open the question of the exact shape of the density for all $\a > 1.$ Observe that the limiting case $\a = 2$ is rather peculiar since it can be elementally shown that its even derivatives never vanish, whereas its odd derivatives vanish only once and at zero. But since the BS property is not closed under pointwise limits, it might be true that $f_{\a,\rho}$ is BS whenever its support is $\rl.$ On the other hand, in spite of Theorem \ref{shapes} (c) we think that for $\a\in (1,2)$ the visually whale-shaped density $f_{\a, 1/\a},$ whose support is a negative half-line, is not in ${\rm WS}_-.$ Indeed, we will see in Proposition \ref{G2} that otherwise it would be ID, and we know that this is not true at least for $\a$ close enough to 2.\\

Our four theorems are proved in Section \ref{MainR}. In the last section, we derive further results related to the analysis of the one-sided free stable densities. First, we analyze in more detail the Kanter random variable $\K_\a$, which plays an important role in the proof of all four theorems. The range $\a <1/5$ is particularly investigated, and two conjectures made in \cite{JS} and \cite{BS} are answered in the negative. A curious Airy-type function is displayed in the case $\a =1/5.$ We also derive the full asymptotic expansion of the densities of $\X_\a$, $\X_{\a,1-1/\a}$ and $1-\F$ at the left end of their support, completing the series representation at infinity (1.16) in \cite{HK}. We then provide some explicit finite factorizations of $\X_\a$ and $\K_\a$ with $\a$ rational in terms of the Beta random variable, and an identity in law for random discriminants on the unit circle is briefly discussed. These factorizations motivate a new identity for the Beta-Gamma algebra, which is derived thanks to a formula of Thomae on the generalized hypergeometric function. Stochastic and convex orderings are obtained for certain negative powers of $\X_\a$, where the free Gumbel law and the exceptional free 1-stable law appear naturally at the limit. We show that some generalizations of the semi-circular random variable $\X_{2,1/2}$ provide a family of examples solving the so-called van Dantzig's problem. Finally, we display some striking properties of whale-shaped functions and densities.
 
\section{Proofs of the main results}
\label{MainR}

\subsection{Preliminaries}

The proofs of all four theorems rely on the following result by Haagerup and M\"oller \cite{HM} who, using a general property of the $S-$transform, have computed the fractional moments of $\X_\a.$ They obtain
\begin{eqnarray*}
\esp [\X_\a^s] & = & \frac{\Gamma(1-s/\a)}{\Gamma (2- (1/\a-1)s)\Gamma (1-s)}\\
& = & \lpa\frac{1}{1+ (1-1/\a)s)}\rpa\times\lpa\frac{\Gamma(1-s/\a)}{\Gamma (1-s) \Gamma (1- (1/\a-1)s)}\rpa
\end{eqnarray*}
for $s <\a.$ Identifying the two factors, we get the following multiplicative identity in law
\begin{equation}
\label{Kant1}
\X_\a\; \elaw \; \U^{1-1/\a}\, \times\, \K_\a,
\end{equation}
where $\U$ is uniform on $(0,1)$ and $\K_\a$ is the so-called Kanter random variable. The latter appears in the following factorization due to Kanter - see Corollary 4.1 in \cite{K}:
\begin{equation}
\label{Kant2}
\Z_\a \; \elaw\; \L^{1-1/\a}\, \times\, \K_\a,
\end{equation}
where $\L$ has unit exponential distribution and $\Z_\a$ is a classical positive $\a-$stable random variable with Laplace transform
$\esp [e^{-\lambda \Z_\a}] = e^{-\lambda^\a}$ and fractional moments
$$\esp[\Z_\a^s]\; =\;\frac{\Gamma(1-s/\a)}{\Gamma (1-s)}$$
for $s <\a.$ Observe that the random variable $\K_\a$ has fractional moments
\begin{equation}
\label{FracSta}
\esp[\K_\a^s]\; =\;\frac{\Gamma(1-s/\a)}{\Gamma(1-(1/\a -1)s)\Gamma (1-s)}
\end{equation}
for $s <\a,$ and in particular a support $[a_\a^{-1}, +\infty)$ which is bounded away from zero, with
$$a_\a\; =\; \a^{-1}(1-\a)^{1-\frac{1}{\a}}\;=\;\lim_{n\to+\infty} \esp[\K_\a^{-n}]^{1/n},$$ 
by Stirling's formula. The density of $\K_\a$ is explicit for $\a = 1/2,$ with 
$$\K_{\frac{1}{2}}\;\elaw\; \frac{1}{4\cos^2(\pi\U/2)}\;\elaw\; \frac{1}{4 \B_{\frac{1}{2},\frac{1}{2}}}$$
and where, here and throughout, $\B_{a,b}$ stands for a standard $\beta(a,b)$ random variable with density
$$\frac{\Ga(a+b)}{\Ga(a)\Ga(b)}\, x^{a-1} (1-x)^{b-1}$$
on $(0,1).$ Plugging this in (\ref{Kant1}) yields easily 
$$\X_{\frac{1}{2}}\;\elaw\;\frac{1}{4 \B_{\frac{1}{2},\frac{3}{2}}}$$
and we retrieve the aforementioned closed expression of $f_{1/2}.$ Several analytical properties of the density of $\K_\a$ have been obtained in \cite{JS, TSEJP}. In particular, Corollary 3.2 in \cite{JS} shows that it is CM, a fact which we will use repeatedly in the sequel.

\begin{remark}
\label{KantFuss}
{\em (a) Specifying Haagerup and M\"oller's result to the negative integers yields 
$$\esp[\X_\a^{-n}]\; =\; \frac{1}{n\a^{-1} +1}\,\binom{n\a^{-1} +1}{n},\quad n\ge 0.$$
The latter is a so-called Fuss-Catalan sequence, and it falls within the scope of more general positive-definite sequences studied in \cite{Mlo10,MPZ}. With the notations of these papers, one has $\X_\a\elaw W_{1/\a,1}^{-1}.$ This implies that $f_\a$ can be written explicitly, albeit in complicated form, for $\a=1/3$ and $\a =2/3$ - see (40) and (41) in \cite{MPZ}. It is also interesting to mention that $\X_{\frac{1}{2}}^{-1}$ has Marchenko-Pastur (or free Poisson) distribution, with density 
$$\frac{1}{2\pi}\sqrt{\frac{4-x}{x}}$$
on $(0,4].$ More generally, Proposition A.4.3 in \cite{BerPat} - see also (8) in \cite{MPZ} - shows that $\X_{\frac{1}{n}}^{-1}$ is distributed  for each $n\ge 2$ as the $(n-1)$-th free multiplicative convolution power of the Marchenko-Pastur distribution.

\vspace{2mm}

(b) The negative integer moments of $\K_\a$ are given by the simple binomial formula
$$\esp[\K_\a^{-n}]\; =\; \binom{n\a^{-1}}{n},\quad n\ge 0.$$
This shows that the law of $\K_\a^{-1}$ is of the type studied in \cite{MP}, more precisely it is $\nu(1/\a,0)$ with the notations therein. By Gauss's multiplication formula - see e.g. Theorem 1.5.2 in \cite{AAR} - and Mellin inversion, this also implies the identity
$$\K_{\frac{1}{3}}^{-1}\;\elaw\;\K_{\frac{2}{3}}^{-2}\;\elaw\; 27\,\B_{\frac{1}{3},\frac{2}{3}}\,(1-\B_{\frac{1}{3},\frac{2}{3}})$$  
in terms of a single random variable $\B_{\frac{1}{3},\frac{2}{3}}.$ In particular, the density of $\K_\a$ can be written in closed form for $\a =1/3$ and $\a =2/3$ as a two-to-one transform of the density of $\B_{\frac{1}{3},\frac{2}{3}}$ - see also Theorems 5.1 and 5.2 in \cite{MP}. As seen above, $\K_{\frac{1}{2}}^{-1}\elaw 4\B_{\frac{1}{2},\frac{1}{2}}$ is arc-sine distributed, with density
$$\frac{1}{\pi\sqrt{x(4-x)}}$$
on $(0,4].$ It is well-known that this is the distribution of the rescaled free independent sum of two Bernoulli random variables with parameter 1/2. It turns out that in general, $\K_{\frac{1}{n}}^{-1}$ is distributed  for each $n\ge 2$ as the $(n-1)$-th free multiplicative convolution power of a free Bernoulli process at time $n/(n-1)$ - see (6.9) in \cite{MP}. 

\vspace{2mm}

(c) The random variable $\K_\a$ can be expressed as the following explicit deterministic transformation of a single uniform variable $\U$ on $(0,1):$ 
\begin{equation}
\label{V}
\K_\alpha\; \elaw\; \frac{\sin(\pi \alpha \U) \sin^{\frac{1-\alpha}{\alpha}}(\pi(1-\alpha) \U)}{\sin^{\frac{1}{\alpha}}(\pi \U)}\cdot 
\end{equation}
This is Kanter's original observation - see Section 4 in \cite{K}, and it will play an important role in the proof of Theorem 3. Notice that the deterministic transformation involved in (\ref{V}) appears in the implicit expression of the densities $f_\a,$ which is given in the second part of Proposition A.I.4 in \cite{BerPat} - see also (11) in \cite{MPZ} for the case when $\a$ is the reciprocal of an integer. There does not seem to exist any computational explanation of this fact. We refer to equation (1) in \cite{D}, and also to Propositions 1 and 2 therein for further results on this transformation.
}
\end{remark}

\subsection{Proof of Theorem 1}

\subsubsection{The case $\a\le 1$} We begin with the one-sided situation $\rho = 1.$ Setting $b_\a = a_\a^{-1}$ we deduce from (\ref{Kant1}) and the multiplicative convolution formula that, for any $x > 0,$
\begin{eqnarray*}
 f_{\a}(x+b_\a) & = & \frac{\a}{1-\a} \int_1^{1+\frac{x}{b_\a}} y^{-\frac{2-\a}{1-\a}} f_{\K_\a}(y^{-1} (x+b_\a))\,dy  \\
         &=& \frac{\a b_\a^{\frac{1}{1-\a}}}{1-\a}\int_0^1 \frac{x}{(b_\a +tx)^{\frac{2-\a}{1-\a}}}\; f_{\K_\a}\!\!\lpa b_\a + \frac{b_\a(1-t)x}{b_\a + tx}\rpa dt.
\end{eqnarray*}
On the one hand, for every $t\in (0,1)$, the function
$$x\; \mapsto\; \frac{(1-t)x}{b_\a + tx}$$
is a Bernstein function - see \cite{SSV}. On the other hand, by the aforementioned Corollary 3.2 in \cite{JS}, the function $z \mapsto f_{\K_\a}(b_\a+z)$ is CM. Hence, by e.g. Theorem 3.7 in \cite{SSV}, the function
$$x\; \mapsto\;  f_{\K_\a}\!\!\lpa b_\a + \frac{b_\a(1-t)x}{b_\a + tx}\rpa$$
is CM, and so is
$$x\; \mapsto\; (b_\a + tx)^{-\frac{2-\a}{1-\a}} f_{\K_\a}\!\!\lpa b_\a + \frac{b_\a(1-t)x}{b_\a + tx}\rpa$$
as the product of two CM functions. Integrating in $t$ shows that $x\mapsto x^{-1}f_{\a}(x+b_\a)$ is CM on $(0,\infty)$ and it is easy to see from Bernstein's theorem that this implies the independent factorization 
$$\X_\a\;\elaw\; b_\a\; +\; \G_2 \,\times\, \Y_\a$$
for some positive random variable $\Y_\a$ where, here and throughout, $\G_t$ stands for a standard $\Gamma(t)$ random variable with density 
$$\frac{x^{t-1}e^{-x}}{\Ga(t)}$$ 
on $(0,+\infty).$ By Kristiansen's theorem \cite{K94}, this shows that $\X_{\alpha}$ is ID.  

To handle the two-sided situation $\rho \in (0, 1) $, we appeal to the following identity in law which was observed in \cite{HK} - see (2.8) therein:
\begin{equation}
\label{HK1}
\X_{\alpha,\rho}\; \elaw \;\X_{1,\rho}\, \times\, \X_\a.
\end{equation}
Since $\X_{1,\rho}$ has a drifted Cauchy law and since the underlying Cauchy process $\{\X_t^{(1, \rho)}\!,\, t \geq 0 \}$ is self-similar with index one, the latter identity transforms into
\begin{equation}
\label{HK2}
\X_{\a, \rho} \; \elaw\; \X_{\X_\a}^{(1, \rho)}
\end{equation}
which is a Bochner's subordination identity. By e.g. Theorem 30.1 in \cite{Sato}, this finally shows that $\X_{\a, \rho}$ is ID for every $\a\in (0,1]$ and $\rho \in [0,1].$

\qed

\begin{remark}
\label{IDone}
{\em (a) The above proof shows that
$$s\;\mapsto\; \frac{\esp[(\X_\a - b_\a)^s]}{\Ga(2+s)}$$
is the Mellin transform of some positive random variable. On the other hand, it seems difficult to find a closed formula for the Mellin transform $\esp[(\X_\a - b_\a)^s],$ except in the case $\a =1/2$ where
$$\esp[(\X_{\frac{1}{2}} - b_{\frac{1}{2}})^s]\; =\; \frac{2^{1-2s}}{\pi}\,\Ga(3/2+s)\,\Ga(1/2 -s),\qquad s\in(-3/2,1/2).$$
When $\a$ is the reciprocal of an integer, there is an expression in terms of the terminating value of a generalized hypergeometric function - see Remark \ref{BG4bis} (c), but we are not sure whether this always transforms into a ratio of products of Gamma functions, as is the case for $\X_\a.$ 

\vspace{2mm}

(b) We believe that $\X_\a - b_\a$ is a $\Ga_{3/2}-$mixture for every $\a\in(0,1),$ that is
$$s\;\mapsto\; \frac{\esp[(\X_\a - b_\a)^s]}{\Ga(3/2+s)}$$
is the Mellin transform of some positive random variable. This more stringent property is actually true for $\a\le 3/4,$ as a consequence of the above proof and Theorem \ref{GGCID} - see Remark \ref{Asympbi} (b).}
\end{remark}

\subsubsection{The case $\a > 1$ and $\rho =1/2$} We first derive a closed expression for the Fourier transform of $\X_{\a,\rho}$, which has independent interest. It was already obtained as Theorem 1.8 in \cite{HK} in a slightly different manner. Our proof is much simpler and so we include it here. Introduce the so-called Wright function
$$\phi(a,b, z)\; =\; \sum_{n\ge 0} \frac{z^n}{\Ga(b+ an)\, n!}$$
with $a > -1, b \in\rl$ and $z\in\CC.$ This function was thoroughly studied in the original articles \cite{W0,W1,W} for various purposes, and is referenced in Formula 18.1(27) in the encyclopedia \cite{EMOT}. It will play a role in other parts of the present paper.

\begin{lemma}
\label{WriFou}
One has
$$\esp[e^{\i t \X_{\a,\rho}}]\; =\; \phi(\a-1,2, -(\i t)^\a e^{-\i\pi\a\rho\, {\rm sgn} (t)}),\qquad t\in\rl.$$
\end{lemma}

\proof The case $\a =1$ is an easy and classic computation, since $\X_{1,\rho}$ has a drifted Cauchy distribution and $\phi(0,2,z) = e^z.$ When $\a\neq 1,$ we first observe that since $\X_{\a,\rho}\elaw- \X_{\a,1-\rho},$ it is enough to consider the case $t >0.$ Combining e.g. Theorem 14.19 in \cite{Sato} and Corollary 1.5 in \cite{HK} yields
$$e^{ -(\i x)^\a e^{-\i\pi\a\rho}}\; =\;\int_0^\infty t\, e^{-t}\,\esp[e^{\i x t^{1-1/\a}\X_{\a,\rho}}]\, dt \; =\; x^{\frac{2\a}{1-\a}}\int_0^\infty t\, e^{-t x^{\frac{\a}{1-\a}}}\,\esp[e^{\i t^{1-1/\a}\X_{\a,\rho}}]\, dt  $$
for all $x >0.$ On the other hand, a straightforward computation implies
$$x^{\frac{2\a}{1-\a}}\int_0^\infty t\, e^{-t x^{\frac{\a}{1-\a}}}\,\phi(\a-1,2, -(\i t^{1-1/\a})^\a e^{-\i\pi\a\rho})\, dt\; =\;e^{ -(\i x)^\a e^{-\i\pi\a\rho}}, \qquad x >0.$$
The result follows then by uniqueness of the Laplace transform.

\endproof

We can now finish the proof of the case $\a > 1, \rho =1/2,$ where the above lemma reads 
$$\esp[e^{\i t \X_{\a,1/2}}]\; =\; \phi(\a-1,2, -\vert t\vert^\a),\qquad t\in\rl.$$
Applying Theorem 1 in \cite{W0} and some trigonometry, we obtain the asymptotic behaviour
$$\phi(\a-1,2, -t^\a)\; \sim\; \kappa_\a\, t^{-3/2}\, e^{\cos(\pi/\a)\, \a (\a-1)^{1/\a -1}\, t}\, \cos(3\pi/2\a\, +\, \sin(\pi/\a)\,\a(\a-1)^{1/\a -1}\, t)$$
as $t\to +\infty,$ for some $\kappa_\a > 0.$ This implies that $t\mapsto \esp[e^{\i t \X_{\a,1/2}}]$ vanishes (an infinite number of times) on $\rl,$ and hence cannot be the characteristic function of an ID distribution - see e.g. Lemma 7.5 in \cite{Sato}.

\qed

\begin{remark}
\label{IDtwo}
{\em (a) It was recently shown in Theorem 1 of \cite{ABS} that for any $a,\beta >0,$ the function $\phi(a,\beta, -z)$ has only positive zeroes on $\CC.$ Combined with Lemma \ref{WriFou}, this entails that the function $t\mapsto  \esp[e^{\i t \X_{\a,\rho}}]$ never vanishes on $\rl$ for $\a > 1$ and $\rho\neq 1/2,$ so that the above simple argument cannot be applied. Nevertheless, we conjecture that $\X_{\a,\rho}$ is not ID for all $\a > 1$ and $\rho\in[1-1/\a,1/\a].$\\

(b) When $\rho = 1/\a,$ Lemma \ref{WriFou} also gives the moment generating function
$$\esp[e^{\lambda \X_{\a,1/\a}}]\; =\; \phi(\a-1,2, \lambda^\a)\; =\; \prod_{n\ge 1} \lpa 1+ \frac{\lambda^\a}{\lambda_{\a,n}}\rpa,\qquad \lambda\ge 0,$$
where $0 < \lambda_{\a,1} < \lambda_{\a,2}\ldots\;$ are the positive zeroes of $\phi(\a-1,2,-z).$ Above, the product representation is a consequence of the Hadamard factorization for the entire function $\phi(\a- 1,2,z)$ which is of order $< 1$ - see again Theorem 1 in \cite{W0}, whereas the simplicity of the zeroes follows from the Laguerre theorem on the separation of zeroes for $\phi(\a- 1,2,z),$ which has genus 0. 

Consider now the random variable
$$\Y_\a\; =\; b_{1/\a}^{-1/\a}\; -\; \X_{\a,1/\a}\; =\; \a(\a-1)^{1/\a -1}\; -\; \X_{\a,1/\a},$$
whose support is $(0,\infty)$ by Proposition A.1.2 in \cite{BerPat}, and whose infinite divisibility amounts to that of $\X_{\a,1/\a}.$ Its log-Laplace transform reads
\begin{eqnarray}
\label{SubZa}
-\log\esp[e^{-\lbd \Y_\a}] & = & \a(\a-1)^{1/\a -1}\lambda\; - \; \sum_{n\ge 1} \log \lpa 1+ \frac{\lambda^\a}{\lambda_{\a,n}}\rpa\\
& = & \int_0^\infty (1- e^{-\lambda^\a x}) \lpa \frac{(\a-1)^{1/\a -1} x^{-1/\a}}{\Ga(1-1/\a)}\, -\, \sum_{n\ge 1} e^{-\lbd_{\a,n} x}\rpa \frac{dx}{x}\nonumber
\end{eqnarray}
where in the second equality we have used Frullani's identity repeatedly and the well-known formula (1) p.viii in \cite{SSV}. Putting everything together shows that $\X_{\a,1/\a}$ is ID if and only if the function on the right-hand side is Bernstein. Unfortunately, this property seems difficult to check at first sight. Observe by Corollary 3.7 (iii) in \cite{SSV} that this function is not Bernstein if the function
$$x\; \mapsto\; \frac{(\a-1)^{1/\a -1} x^{-1/\a}}{\Ga(1-1/\a)}\, -\, \sum_{n\ge 1} e^{-\lbd_{\a,n} x}$$
takes negative values on $(0,\infty)$, but this property seems also difficult to study. A lengthy asymptotic analysis which will not be included here, shows that it converges at zero to some positive constant. \\

(c) Rewriting equation (\ref{SubZa}) as
$$\lbd^{1/\a}\; = \; b_{1/\a}\lpa  \sum_{n\ge 1} \log \lpa 1+ \frac{\lambda}{\lambda_{\a,n}}\rpa\; -\;\log\esp[e^{-\lbd^{1/\a} \Y_\a}]\rpa,$$
we obtain the factorization
$$\Z_{1/\a}\;\elaw \; b_{1/\a}\lpa  \sum_{n\ge 1}\; {\rm Exp} (\lambda_{\a,n})\; +\; \Z^{(1/\a)}_{\Y_\a}\rpa$$
where $\{\Z^{(1/\a)}_t, \, t\ge 0\}$ is the $(1/\a)-$stable subordinator and all quantities on the right-hand side are independent. This identity is similar to that of the Lemma in \cite{TSPAMS}, except that the parameters $\lbd_{\a,n}$ of the exponential random variables are not explicit.

}
\end{remark}
\subsection{Proof of Theorem 2} 
\subsubsection{The case $\rho =1$} Here, we need to show that the law of $\X_\a$ is a true GGC. To do so, we first observe that by (\ref{Kant1}) and some rearrangements, one has 
\begin{equation}
\label{survive}
\esp[e^{-\lbd \X_\a}]\; =\; \frac{\a\, \lbd^{\frac{\a}{1-\a}}}{1-\a}\,\int_\lbd^\infty \esp[e^{-x \K_\a}]\, x^{\frac{1}{\a-1}}\, dx,\qquad \lbd \ge 0.
\end{equation}
A combination of Theorem 6.1.1 and Properties (iv) and (xi) p.68 in \cite{B} imply then that it is enough to show that the law of $\K_\a$ itself is a GGC. Alternatively, one can use the main result of \cite{BJTP}, since it is easily seen that $\U^{1-1/\a}$ has a GGC distribution. 
To analyze the law of $\K_\a,$ we use the identity in law
\begin{equation}
\label{Kant3}
\K_\a\;\elaw\; \K_{1-\a}^{\frac{1}{\a} -1},
\end{equation}
a consequence of (\ref{FracSta}) which shows that both random variables have the same fractional moments. Plugging (\ref{Kant3}) again into (\ref{Kant2}) implies that the Laplace transform of $\K_{1-\a}$ is the survival function of the power transformation $\Z_\a^{\frac{-\a}{1-\a}}.$ In other words, one has
\begin{equation}
\label{LapSur}
\esp[e^{-x\K_{1-\alpha}}]\; =\; \pb[\L \geq x \K_{1-\alpha}]\; =\; \pb[\Z_\a^{-\frac{\alpha}{1-\alpha}} \geq x],\qquad x \ge 0.
\end{equation}
Setting $F_\a(x)$ for the function defined in (\ref{LapSur}), we next observe that since $\K_\a$ has a CM density and support $[b_\a, +\infty),$ this function $F_\a$ has by Theorem 9.5 in \cite{SSV} an analytic extension on $\CC\setminus (-\infty,0]$ which is given by
\begin{equation}
\label{RepFa}
F_\alpha(z)\; =\; \exp -\lcr b_{1-\a} z \, +\, \int_0^\infty \frac{z}{t+z}\,\frac{\theta_\a(t)}{t}\, dt\rcr
\end{equation}
for some measurable function $\theta_\a: (0,\infty)\to [0,1]$ such that $\int_0^1 \theta_\alpha(t) t^{-1}\,dt <\infty.$ See also Theorem 51.12 in \cite{Sato}. Applying now Theorem 8.2 (v) in \cite{SSV}, we see that the GGC property of $\K_{1-\a}$ is equivalent to the non-decreasing character of $\theta_\a$ on $(0,\infty),$ and the following proposition allows us to conclude the proof of the case $\rho = 1.$

\begin{proposition}
\label{ThInc}
The function $\theta_\a$ has a continuous version on $(0,\infty)$, which is non-decreasing for every $\a\in[1/4,1).$
\end{proposition}

\proof

The analysis of $\theta_\a$ depends, classically, on the behaviour of $F_\a$ near the cut. Assume for a moment that $\theta_\a$ is continuous. For every $r >0$ and $\delta \in (0,1),$ we have, after some simple rearrangements,
\begin{eqnarray*}
\frac{F_\a( r e^{\i\pi(1-\delta)} )}{F_\a(r e^{\i\pi(\delta-1)})} & = & \exp - 2\i \lcr\sin(\pi\delta) b_{1-\a} r \, +\,\int_0^\infty \frac{\sin(\pi(1-\delta))\, \theta_\a(rt)}{1+ 2\cos(\pi(1-\delta))t + t^2}\, dt\rcr\\
& = & \exp - \lcr 2\i \sin(\pi\delta) b_{1-\a} r \, +\,2\i\pi \esp \lcr\theta_\a(r (\X_{1, 1-\delta})^+)\rcr\rcr\\ 
&\to & \exp -2\i\pi\theta_\a(r)
\end{eqnarray*}
as $\delta\to 0$ since $\X_{1,1-\delta}\to\Un$ in law as $\delta\to 0$ and $\theta_\a$ is bounded continuous. On the other hand, it follows from the third expression of $F_\a$ in (\ref{LapSur}) and the first formula of Corollary 1 p.71 in \cite{Z}, after a change of variable, that
$$F_\a(x)\;=\; 1+ \frac{1}{2\i\pi} \int_0^\infty e^{-t}\lpa e^{-e^{\i\pi \alpha} t^\alpha x^{1-\alpha}} - e^{-e^{-\i\pi \alpha} t^\alpha x^{1-\alpha}}\rpa\frac{dt}{t} , \qquad x>0.$$
The analytic continuations of $F_\a$ near the cut are then expressed, changing the variable backwards, as
$$F_\a(r e^{\i \pi}) \;=\; 1 + \frac{1}{2\i\pi} \int_0^\infty e^{-ru}\lpa e^{r u^\alpha }- e^{r u^\alpha e^{-2\i\pi\alpha}}\rpa \frac{du}{u}$$
and
$$F_\a(r e^{-\i \pi})\; =\; 1 + \frac{1}{2\i\pi} \int_0^\infty e^{-ru}\lpa e^{r u^\alpha e^{2\i\pi\alpha}}- e^{r u^\alpha }\rpa\frac{du}{u}\; = \; \overline{F_\a (r e^{\i\pi})}.$$
Therefore, we obtain
$$\frac{F_\alpha(re^{\i\pi})}{F_\alpha(re^{-\i\pi})}\; =\; e^{-2\i\pi \eta_\alpha(r)}$$
for every $r > 0,$ with the notation
$$\eta_\a (r)\; = \;  \frac{1}{\pi}\arg [F_\a(re^{-\i\pi})].$$
Since
$$\Im (F_\a(re^{-\i\pi}))\; = \; \frac{1}{2\pi} \int_0^\infty e^{-ru}\,\Re\lpa e^{ru^\a} - e^{r u^\a e^{2\i\pi\a}} \rpa \frac{du}{u}\; >\; 0$$ 
for every $r >0,$ the function $\eta_\a$ takes its values in $[0,1]$ and is clearly continuous. By construction, the functions $t^{-1}\eta_\a(t)$ and $t^{-1}\theta_\a(t)$ have the same Stieltjes transform, and it follows by uniqueness that $\theta_\a$ has a continuous version, which is $\eta_\a.$\\

It remains to study the monotonous character of $\eta_\a$ on $(0,\infty).$ A first observation is that, expanding the exponentials inside the brackets and using the complement formula for the Gamma function, the following absolutely convergent series representation holds:
\begin{equation}
\label{Wrgt}
F_\a (re^{-\i\pi})\; =\; \sum_{n\ge 0} \frac{z^ne^{\i\pi n\a}}{n!\Gamma(1-n \a)}\; =\; \phi (-\a, 1, z e^{\i\pi\a})
\end{equation}
with $z = r^{1-\a}.$ In particular, the function
$$r\;\mapsto\;r^{\a-1}\,\Im\lpa F_\a (re^{-\i\pi})\rpa\; =\;\frac{1}{\pi} \sum_{n\ge 1} \frac{\Gamma(n\a)}{n!}\,r^{(n-1)(1-\a)}$$
is absolutely monotonous on $(0,\infty),$ and the non-decreasing character of $\theta_\a$ will hence be established as soon as $r\mapsto r^{\a-1} \Re(F_\a(re^{-\i\pi}))$ is non-increasing on $(0,\infty).$ We use the representation
$$\Re(F_\a(r e^{-\i\pi}))\;=\;1\;+\; \frac{1}{2\pi}\int_0^\infty e^{-r^{1-1/\alpha}u}\,\Im ( e^{ u^\alpha e^{2\i\pi\alpha}})\,\frac{du}{u}$$
and divide this last part of the proof into three parts.

\begin{itemize}

\item The case $\alpha\in[1/2,1).$ If $\a =1/2,$ we simply have $\Re(F_{1/2}(re^{-\i\pi})) \equiv 1.$ If $\a > 1/2$ we rewrite, using again the first part of Corollary 1 p.71 in \cite{Z},
\begin{eqnarray*}
\Re(F_\a(r e^{-\i\pi})) & = &\frac{1}{2}\; +\; \frac{1}{2}\lpa 1-\frac{1}{\pi}\int_0^\infty e^{-r^{1-1/\a}u}\,\Im(e^{-u^\alpha e^{-i\pi\rho\alpha}})\,\frac{du}{u}\rpa \\
& = & \frac{1}{2}\lpa 1\, +\, \pb[\Z_{\alpha,\rho}\leq r^{1-1/\a}]\rpa
\end{eqnarray*}
where $\rho=2-1/\a \in (0,1)$ and $\Z_{\a,\rho}$ is as in Lemma \ref{WriFou}  a real $\a-$stable random variable with positivity parameter $\rho.$
Thus, $\Re(F_\a(r e^{-\i\pi}))$ decreases (from 1 to $1/2\a$) on $(0,\infty)$ and $r^{\a-1} \Re(F_\a(r e^{-\i\pi}))$ also decreases on $(0,\infty),$ as required.\\

\item The case $\a\in[1/3,1/2).$ Setting $\rho = 1/\a -2\in (0,1]$ and using the same notation as in the previous case, we rewrite
\begin{eqnarray*}
\Re(F_\a(re^{-\i\pi})) & = & 1 \; +\; \frac{1}{2\pi}\int_0^\infty  e^{-r^{1-1/\a}u}\,\Im(e^{-u^\alpha e^{-\i\pi\rho\a}})\,\frac{du}{u} \\
                          &=& 1\; +\; \frac{1}{2}\,\pb[\Z_{\a,\rho} \ge r^{1-1/\a}] \\ & = & 1\;+\;\frac{\rho}{2}\,\pb[\W_{\a,\rho}^{-\a} \leq r^{1-\a}]
\end{eqnarray*}
where $\W_{\a,\rho}\elaw\Z_{\a,\rho}\,\vert\,\Z_{\a,\rho} > 0$ is the cut-off random variable defined in Chapter 3 of \cite{Z}. Observe that here, the function $r\mapsto\Re(F_\a(re^{-\i\pi}))$ increases. Setting $h_{\a,\rho}$ for the density function of $\W_{\a,\rho}^{-\a}$ on $(0,\infty),$ we get after a change of variable
$$r^{\alpha -1}\,\Re(F_\a(r e^{-\i\pi})) \; =\; r^{\alpha -1}\;+\; \frac{\rho}{2}\,\int_0^1 h_{\alpha,\rho}(r^{1-\alpha} x)\, dx,$$
and it is hence sufficient to prove that the function $h_{\a,\rho}$ is non-increasing on $(0,\infty).$ 
Using the expression for the Mellin transform of $\W_{\a,\rho}$ given at the bottom of p.186 in \cite{Z} together with the complement and multiplication formul\ae\, for the Gamma function, we obtain
\begin{eqnarray*}
\esp[\W_{\a,\rho}^{-\a s}] & = & \frac{\Gamma(1+ s)}{\Gamma(1+\a\rho s)}\,\times\, \frac{\Gamma(1-\a s)}{\Gamma(1-\a\rho s)} \\
& = &\frac{2^s}{1+s}\,\times\, \frac{\Gamma(3/2+s/2)}{\Gamma(3/2)}\,\times\, \frac{\Gamma(1+s/2)}{\Gamma(1+\a\rho s)}\,\times\,  \frac{\Gamma(1-\a s)}{\Gamma(1-\a\rho s)}
\end{eqnarray*}
for every $s\in (-1, 1/\a).$ Identifying the factors and using $\a\rho < 1/2,$ this implies the identity in law
$$\W_{\a,\rho}^{-\a}\;\elaw \; 2\U\,\times\, \sqrt{\G_{3/2}}\,\times\,\lpa\frac{\Z_{\rho}}{\Z_{2\alpha\rho}}\rpa^{\a\rho}$$
where all factors on the right hand side are assumed independent. Hence, $\W_{\a,\rho}^{-\a}$ admits $\U$ as a multiplicative factor and by Khintchine's theorem, its density is non-increasing on $(0,\infty).$  \\

\item The case $\a\in [1/4, 1/3).$ Contrary to the above, the argument is here entirely analytic. We consider
 \begin{eqnarray*}
G_\a(r)\; =\; \Re(F_\a(r^{\frac{\a}{\a-1}}e^{-\i\pi})) & =& 1\;+\; \frac{1}{2\pi}\int_0^\infty e^{-ru}\,\Im(e^{ u^\alpha e^{2\i\pi\alpha}})\,\frac{du}{u}\\
                                                       &=& 1\; +\; \frac{1}{2\pi}\int_0^\infty e^{\cos(2\pi\alpha)u^\a-ru} \sin(\sin(2\pi \alpha)u^\a)\,\frac{du}{u} \\
&=& 1\; +\; \frac{1}{2\pi\a}\int_0^\infty g_{\a,r}(t)\,\sin(t)\, dt
\end{eqnarray*}
where
$$t\;\mapsto\; g_{\a,r}(t) \;=\; t^{-1}e^{\cot(2\pi\alpha) t -r(\sin(2\pi\alpha))^{-1/\alpha} t^{1/\a}}$$ 
decreases on $(0, +\infty).$ For every $k\ge 0$ we have
$$\int_{2k\pi}^{2(k+1)\pi} g_{\a, r}(t)\, \sin (t)\, dt \; = \;\int_0^{\pi} (g(t+ 2k\pi) - g(t+ (2k+1)\pi)) \,\sin(t)\, dt \;>\; 0,$$
so that $G_{\alpha}(r) > 1$ for every $r >0.$ We next compute
\begin{eqnarray*}
(r^\a G_\a(r))' &=& \alpha r^{\alpha -1}(G_{\alpha}(r)\, +\, \frac{r}{\alpha} G'_{\alpha}(r)) \\
      &>&  \alpha r^{\alpha -1}\lpa 1 - \frac{r}{2\pi\alpha}\int_0^\infty e^{\cos(2\pi\alpha)u^\a - ru} \sin(\sin(2\pi \alpha)u^\a)du\rpa   \\
      &=& \frac{r^\a}{2\pi} \int_0^\infty e^{-ru} \lpa 2\pi \alpha - e^{\cos(2\pi\alpha) u^\a} \sin(\sin(2\pi \alpha) u^\a)\rpa du \; >\; 0, 
\end{eqnarray*}
since $2\pi \alpha > 1 \geq e^{\cos(2\pi\alpha) u^\a} \sin(\sin(2\pi \alpha)u^\a)$ for every $u> 0.$ Changing the variable backwards, this finally shows that $r\mapsto r^{\alpha - 1}\Re(F_\a(r e^{-\i\pi}))$ decreases on $(0,\infty).$ 

\qed

\end{itemize}

\begin{remark}
\label{Rho1}
{\em (a) The above argument shows that the survival function $x\mapsto\pb[\Z_\a^{-\frac{\alpha}{1-\alpha}} \geq x]$ is HCM for every $\a\ge 1/4$, with the terminology of \cite{B}. A consequence of Corollary \ref{NegKa} is that this is not true anymore for $\a < 1/5,$ and we believe  - see Conjecture \ref{KaGGC} - that the right domain of validity of this property is $\a\in [1/5,1)$. The more stringent property that $\Z_\a^{-\frac{\alpha}{1-\alpha}}$ is a HCM random variable for $\a\le 1/2$ was conjectured in \cite{PB} and some partial results were obtained in \cite{PB,BS}. In \cite{F}, it is claimed that this latter property holds true if and only if $\a\in[1/3,1/2].$ 

\vspace{2mm}

(b) The analytical proof for the case $\a\in[1/4,1/3)$ conveys to the case $\a\in[1/3,1/2).$ Nevertheless, it is informative to mention the probabilistic interpretation of $\Re(F_\a(re^{-\i\pi}))$ for $\a\in[1/3,1/2).$ Simulations show that this function oscillates for $\a < 1/3.$ See also Section 4.2 for a striking similarity between the cases $\a =1/3$ and $\a =1/5.$

\vspace{2mm}

(c) We do not know if the representation (\ref{RepFa}) holds for the Laplace transform of $\X_\a.$ Since the latter is a $\Ga_2-$mixture we obtain, similarly as above,
$$\esp[e^{-x \X_\a}] \; =\;e^{-b_{1-\a} x} \int_0^{\infty}\frac{\nu_\a(dt)}{(x+t)^2}$$
for some positive measure $\nu_\a$ on $[0,+\infty).$ This representation would suffice if we could show that the generalized Stieltjes functions on the right-hand side is the product of two standard Stieltjes functions, applying Theorem 6.17 in \cite{SSV} as in the proof of Theorem 9.5 therein. However, this is not true in general, for example when $\nu_\a$ is the sum of two Dirac masses. Observe that in the other direction, the product of two Stieltjes functions is a generalized Stieltjes function of order 2 - see Theorem 7 in \cite{KP}. With the notation of \cite{KP}, we believe that the exact Stieltjes order of $\esp[e^{-x \X_\a}]$ is actually $3/2,$ which however does not seem of any particular help for (\ref{RepFa}). Alternatively, because of (\ref{survive}) one would like to prove that if $f$ has representation (\ref{RepFa}), then so has $x\mapsto \int_x^\infty f(y) dy$. This is true in the GGC case by Property xi) p.68 in \cite{B}, but we were not able to prove this in general.}
\end{remark}

\subsubsection{The case $\rho < 1$} The case $\rho = 0$ follows from $\X_{\a,0} \elaw -\X_\a.$ For $\rho\in (0,1)$ we appeal to (\ref{HK2}), the previous case, and the Huff-Zolotarev subordination formula which is given e.g. in Theorem 30.1 of \cite{Sato}. Since the law of $\X_\a$ is a GGC for $\a\le 3/4$, its Laplace transform reads
$$\esp[e^{-\lbd\X_\a}]\; =\;\exp-\lcr b_\a\lbd\; +\; \int_0^\infty (1-e^{-\lbd x})\, k_\a(x)\,\frac{dx}{x}\rcr$$ 
for some CM function $k_\a.$  Formula (30.8) in \cite{Sato} and the closed expression of the density of $\X_{1,\rho}$ imply that the L\'{e}vy measure $\nu_{\a,\rho}$ of $\X_{\a,\rho}$ has density
\begin{eqnarray*}
\psi_{\a,\rho}(x) & = & b_\a\psi_{1,\rho} (x)\;+\; \frac{\sin (\pi\rho)}{\pi} \int_0^\infty \frac{k_\a (u)}{x^2 + 2\cos (\pi\rho) xu  + u^2}\, du\\
& = &  \frac{\sin (\pi\rho)}{\pi \vert x \vert} \lpa \frac{b_\a}{\vert x\vert}\;+\;  \int_0^\infty \frac{k_\a (\vert x \vert u)}{1 + 2\cos (\pi\rho)\,\text{sgn}(x)\, u  + u^2}\, du\rpa
\end{eqnarray*}
over $\rl^*,$ where the closed expression for $\psi_{1,\rho}$ can be deduced e.g. from Theorem 14.10 and Lemma 14.11 in \cite{Sato}. Both functions $x\psi_{\a,\rho}(x)$ and $x\psi_{\a,\rho}(-x$) are hence CM on $(0,\infty).$ 

\qed

\begin{remark}
\label{Triplet}
{\em Since $b_\a > 0$ and the ID random variable $\X_{1,\rho}$ has no Gaussian component, the Huff-Zolotarev subordination formula shows that $\X_{\a,\rho}$ does not have a Gaussian component either, and that for $\rho\in (0,1)$ its L\'evy measure is such that
$$\int_{\vert x\vert \le 1} \vert x\vert\, \nu_{\a,\rho} (dx)\; = \; +\infty.$$
With the terminology of \cite{Sato} - see Definition 11.9 therein, this means that the L\'evy process associated with $\X_{\a,\rho}$ is of type C. This contrasts with the classical $\a-$stable L\'evy process which is of type B for $\a <1.$ When $\rho = 1$ and $\a\le 3/4,$ the GGC property shows that the L\'evy process corresponding to $\X_\a$ is of type B. We believe that this is true for all $\a\in (0,1),$ but this cannot be deduced from the sole $\Ga_2-$mixture property established in Theorem 1.}
\end{remark}

\subsection{Proof of Theorem 3} It is well-known and easy to see from the Voiculescu transform 
$$\phi_{1,1/2} (z)\; =\; -\i$$ 
that the free independent sum of $\X_{1,1/2}$ with any random variable is also a classical independent sum. Hence, the ID character of $\C_{a,b}$ follows from that of $\F,$ which is a consequence of Theorem \ref{freeID} and the convergence in law
\begin{equation}\label{conv T}
\frac{(1-\a)^{1-\a} - \X_\a}{1-\a}\;\claw\;\F\qquad\mbox{as $\a\uparrow 1,$}
\end{equation}
the latter being easily obtained in comparing the two Voiculescu transforms. This concludes the first part of the theorem. Moreover, it is clear that neither $\X_{1,1/2}$ nor $\F$, whose support is a half-line by Proposition A.1.3 in \cite{BerPat}, have a Gaussian component, and this property conveys hence to $\C_{a,b}.$ Finally, since the L\'evy measure of $\X_{1,1/2}$ is
$$\frac{1}{\pi x^2}\,\Un_{\{x\neq 0\}}$$
as seen in the above proof, we are reduced to show by independence and scaling that the L\'evy measure of $\F$ has density
$$\frac{1}{x^2}\lpa 1\, -\, \frac{\vert x\vert\,e^{-2\vert x\vert}}{1- e^{-\vert x\vert}}\rpa\Un_{\{ x < 0\}}.$$
This last computation will be done in two steps. Consider the random variable
$$\W\; =\; \frac{\sin (\pi \U)}{\pi \U}\, e^{\pi \U \cot(\pi \U)}$$
and the exceptional $1$-stable random variable $\S$ characterized by
$$\bE[e^{s \S}]\;=\; s^s, \qquad s > 0.$$
 
\begin{proposition}
\label{prop:SF} 
One has the identities
$$\S\; \elaw\; \log\L\, +\, \log\W\qquad\quad\mbox{and}\qquad\quad\F\; \elaw\; \log\U\, +\, \log\W.$$
\end{proposition}

\proof 
 We begin with the first identity. Using (\ref{Kant2}), we decompose
\begin{equation}
\label{DD}
\frac{(1-\a)^{1 -\a} - \Z_\a}{1-\a}\; \elaw\; \K_\a\,\times\lpa \frac{1- \L^{1-\frac{1}{\a}}}{1-\a}\rpa\,+ \,\lpa\frac{(1-\a)^{1 -\a} - \K_\alpha}{1-\alpha}\rpa.
\end{equation}
On the one hand, a comparison of the two moment generating functions yields 
$$\frac{(1-\a)^{1 -\a} - \Z_\a}{1-\a}\;\claw\;\S\qquad\mbox{as $\a\uparrow 1.$}$$
On the other hand, the right-hand side of (\ref{DD}) is a deterministic transformation, depending on $\a,$ of $(\L,\U)$ independent. It is easy to see from (\ref{V}) that
$$\K_\a\,\times\lpa \frac{1- \L^{1-\frac{1}{\a}}}{1-\a}\rpa\;\pslaw\;\log \L\qquad\mbox{as $\a\uparrow 1.$}$$ 
To study the second term, we use the elementary expansions
\begin{align}
\nonumber
\sin(\pi \alpha \U)\; &=\;  \sin(\pi \U)\, + \,(\alpha-1)\pi\U\cos(\pi \U)\, + \,O((1-\alpha)^3) \\
\nonumber\sin^{\frac{1-\alpha}{\alpha}}(\pi(1-\alpha) \U) \; &=\; 1+(1-\alpha) \log \sin (\pi \U (1-\alpha))\, +\, O((1-\alpha)^2\log^2(1-\alpha)) \\
\nonumber\sin^{\frac{1}{\alpha}}(\pi \U)\; &=\; \sin(\pi \U) (1+(1-\alpha) \log\sin(\pi \U))\, +\, O((1-\alpha)^2) \\
\nonumber(1-\alpha)^{1-\alpha}\; &=\;  1+(1-\alpha)\log(1-\alpha)\, +\, O((1-\alpha)^2 \log^2(1-\alpha))
\end{align}

\noindent
which, combined with (\ref{V}), yield the almost sure asymptotics 
$$\frac{(1-\a)^{1 -\a} - \K_\alpha}{1-\alpha}\; =\;  \log\left(\frac{\sin(\pi \U)}{\pi \U}e^{\pi \U \cot(\pi \U)}\right)\, +\, O((1-\alpha)\log^2(1-\alpha)).$$
Putting everything together completes the proof of the first identity. The second one is derived exactly in the same way, using (\ref{Kant1}) and (\ref{conv T}).

\endproof

\begin{remark}
\label{prop:SFbis}
{\em (a) The first identity in Proposition \ref{prop:SF} is actually the consequence of an integral transformation due to Zolotarev - see (2.2.19) with $\beta =1$ in \cite{Z}. We have offered a separate proof which is perhaps clearer, and which enhances the similarities between the free and the classical case echoing those between (\ref{Kant1}) and (\ref{Kant2}). Observe in particular the identity
\begin{equation}
\label{SF}
\S\; \elaw\; \F\; + \;\log \G_2
\end{equation}
reminiscent of Corollary 1.5 in \cite{HK}, and which is a consequence of Proposition \ref{prop:SF} and the standard identities
\begin{equation}\label{eq:LU}
\L^\beta \; \elaw \;\U^\beta\;\times\;\G_2^\beta
\end{equation}
valid for every $\beta\in\rl^*$ and their limit as $\beta\to 0,$ which is
\begin{equation}\label{eq:LU2}
\log\L\;\elaw\; \log\U\; +\; \log\G_2.
\end{equation}  

\vspace{2mm}

(b) It is interesting to look at these standard identities \eqref{eq:LU} and \eqref{eq:LU2} in the context of extreme value distributions.  Indeed, the three classical extreme distributions are Fr\'echet $\L^\beta$ for $\beta < 0,$ Weibull $-\L^\beta$ for $\beta > 0$ and Gumbel $-\log \L $ for $\beta\to 0$, whereas the free counterparts are $\U^\beta$ for $\beta < 0,$ $-\U^\beta$ for $\beta > 0$ and $-\log \U $ for $\beta\to 0$ according to the classification of \cite{BAV06}.

\vspace{2mm}

(c) Recently Vargas and Voiculescu have introduced Boolean extreme value distributions \cite{VV}. The result is the Dagum distribution, which is indexed by $\beta >0$ and has density function 
$$\frac{x^{1/\beta-1}}{\beta(1+x^{1/\beta})^2}$$
on $(0,\infty).$ Hence, the Dagum distribution is the law of 
$$(\U^{-1} -1)^\beta\;\elaw\;\lpa\frac{\L}{\L}\rpa^\beta$$
which is the independent quotient of two Fr\'echet distributions, and an example of the generalized Beta distribution of the second kind (GB2). On the other hand, by Proposition 4.12 (b) in \cite{AH16a}, the Boolean $\a-$stable distribution has for $\a\le 1$ the law of the independent quotient
$$\frac{\Z_{\alpha,\rho}}{\Z_\alpha}$$
and it is interesting to notice that by Zolotarev's duality - see (3.3.16) in \cite{Z} - and scaling, the positive part of this random variable is distributed as
$$\lpa\frac{\Z_{\alpha\rho}}{\Z_{\alpha\rho}}\rpa^\rho\;\claw\; \lpa\frac{\L}{\L}\rpa^{\frac{1}{\a}}\qquad\mbox{as $\rho\to 0.$}$$
Finding an interpretation about why such quotients appear in those two Boolean cases is left to future work. 

\vspace{2mm}

(d) The second identity in Proposition \ref{prop:SF} can be rewritten as
$$e^{\F}\; \elaw\; \U \;\times\; \W.$$  
In \cite{AH}, it is pointed out that the law of $e^{\F}$ is the Dykema-Haagerup distribution, which appears as the eigenvalue distribution of $A_N^* A_N$ as $N\to\infty$, where $A_N$ is an $N\times N$ upper-triangular random matrix with independent complex Gaussian entries - see \cite{DH04}. 

\vspace{2mm}

(e) It follows from Euler's product and summation formul\ae\, for the sine and the cotangent that $\log\W$ is a decreasing concave deterministic transformation of $\U.$ This implies easily that $\log \W$ has an increasing density on its support which is $(-\infty,1].$ In particular, $\log \W$ is unimodal. Besides, since the densities of $\log\U$ and $\log\L$ are clearly log-concave on the interior of their support, applying Theorem 52.3 in \cite{Sato} we retrieve the known facts that $\S$ and $\F$ are unimodal random variables.}
\end{remark}

Our second step is to compute the Mellin transform of $\W.$
 
\begin{proposition}
\label{MellW}
One has
$$\bE[\W^s]\; =\; \frac{s^s}{\Gamma(1+s)}\; =\; \exp \lcr s \, -\, \int_0^\infty (1-e^{-sx}) \lpa 1 - \frac{x}{e^x-1}\rpa \frac{dx}{x^2}\rcr$$
for all $s > 0.$
\end{proposition}
\proof
The first equality follows from 
$$s^s\; =\; \bE[e^{s \S}]\; =\; \bE[\L^s]\,\bE[\W^s]\; =\; \Ga(1+s)\,\bE[\W^s], \qquad s > 0,$$
a consequence of the first identity in Proposition \ref{prop:SF}. To get the second one, we proceed as in the proof of Lemma 14.11 of \cite{Sato} and start from Frullani's identity
$$\log s\; =\; \int_0^\infty \frac{e^{-x} - e^{-sx}}{x}\,dx$$
which transforms, dividing the integral at 1 and making an integration by parts, into
$$ s\log s\; =\; \int_0^\infty (e^{-sx} - 1 + sx \Un_{\{x \le 1\}})\,\frac{dx}{x^2}\; -\; s\lpa\int_0^\infty (e^{-x} - 1 + x \Un_{\{x \le 1\}})\,\frac{dx}{x^2}\rpa.$$
On the other hand, it is well-known - see e.g. Proposition 4 (a) in \cite{Y} - that
$$\log\Ga(1+s)\; =\; -\gamma s \; +\; \int_0^\infty (e^{-sx} - 1 + sx)\,\frac{dx}{x(e^x -1)}$$
where $\gamma = -\Ga'(1)$ is Euler's constant. Combining the two formul\ae\, yields
\begin{eqnarray*}
\log\esp[\W^s] & = & cs \; +\;  \int_0^\infty (e^{-sx} - 1 + sx \Un_{\{x \le 1\}}) \lpa 1 - \frac{x}{e^x-1}\rpa \frac{dx}{x^2}\\
& = & {\tilde c}s \; -\;  \int_0^\infty (1 - e^{-sx}) \lpa 1 - \frac{x}{e^x-1}\rpa \frac{dx}{x^2}
\end{eqnarray*}
where $c, {\tilde c}$ are two constants to be determined. But it is clear that ${\tilde c}$ is the right end of the support of $\log\W$ which we know, by Remark \ref{prop:SFbis} (c), to be one. Alternatively, one can use Binet's formula
$$\gamma\; =\; \int_0^\infty \lpa \frac{e^{-x}}{1-e^{-x}} \, -\, \frac{e^{-x}}{x}\rpa dx,$$
which is 1.7.2(22) in \cite{EMOT} for $z =1,$ and rearrange the different integrals, to retrieve ${\tilde c} = 1.$ This completes the proof.

\endproof

We can now finish the proof of Theorem \ref{Free1}. Putting together Propositions  \ref{prop:SF} and \ref{MellW}, we get
\begin{eqnarray*}
\log\esp[e^{s\F}] \; = \; \log\esp[\U^s]\; +\; \log\esp[\W^s]
& = & -\log (1+s)\; +\; \log\esp[\W^s]\\
& = & s \; -\; \int_0^\infty (1 - e^{-sx}) \lpa 1 - \frac{x\,e^{-2x}}{1-e^{-x}}\rpa \frac{dx}{x^2}
\end{eqnarray*}
where the third equality follows from rearranging Frullani's identity and the second equality in Proposition \ref{MellW}. All of this shows that the ID random variable $\F$ has support $(-\infty,1]$ - in accordance with Proposition A.1.3 in \cite{BerPat}, and that its L\'evy measure has density
$$\frac{1}{x^2}\lpa 1\, -\, \frac{\vert x\vert\,e^{-2\vert x\vert}}{1- e^{-\vert x\vert}}\rpa\Un_{\{ x < 0\}}$$
as required.

\qed

\begin{remark}
\label{MellWbis}
{\em (a) The first equality in Proposition \ref{MellW} shows that $\W$ has the distribution $\nu_0$ studied in Theorem 6.1 of \cite{Mlo10}. This distribution also appears in Sakuma and Yoshida's limit theorem - see \cite{SY13}. Finally, combining this equality and the second identity in Proposition \ref{prop:SF} implies
$$\esp[e^{s\F}]\; =\; \frac{s^s}{\Gamma(2+s)}$$ 
for all $s >0,$  which was previously obtained in \cite{AH} by other methods, and will be used henceforth.

\vspace{2mm}

(b) It is easy to see that the function
$$x\;\mapsto \; \frac{1}{x}\; -\; \frac{1}{e^x -1}$$
decreases from $1/2$ to zero on $(0,\infty).$ By Corollary 15.11 in \cite{Sato}, this shows that $\log \W$ is SD. A further computation yields
\begin{equation}
\label{CMJ}
\frac{1}{x^2}\lpa 1\; -\; \frac{x}{e^x -1}\rpa\; =\; \int_0^\infty e^{-ux} (u -[u]) \, du, \qquad x > 0.
\end{equation}
This implies that $\log\W$ has CM jumps and that, by Theorem 3, so does $\F$ whose L\'evy measure has density
$$\frac{e^{-\vert x\vert}}{\vert x\vert}\; +\; \frac{1}{x^2}\lpa 1\; -\; \frac{\vert x\vert}{e^{\vert x\vert} -1}\rpa\; =\; \int_0^\infty e^{-u\vert x\vert} (u -[u-1]_+) \, du, \qquad x < 0.$$
By Theorem 51.12 in \cite{Sato}, the latter computation also implies that the law of the positive random variable $1 -\log\W$ is a mixture of exponentials (ME) viz. it has a CM density, which improves on Remark \ref{prop:SFbis} (d) and will be used henceforth. Reasoning as in Corollary 3.2 in \cite{JS} finally implies that the law of
$$\frac{1}{\W}\; -\;\frac{1}{\e}$$
is an ME as well.

\vspace{2mm}

(c) Making an integration by parts in (\ref{CMJ}) yields
$$\frac{1}{x}\lpa 1\; -\; \frac{x}{e^x -1}\rpa\; =\; \int_0^\infty e^{-ux} \lpa du - \sum_{n\ge 1} \delta_n (du) \rpa$$
where $\delta$ stands for the Dirac mass. By (7.1.5) in \cite{B}, this implies that the law of $\log \W$ is not a GGC, and the same is true for $\F$ because
$$\frac{1}{x}\lpa 1\; -\; \frac{x\, e^{-2x}}{1 - e^{-x}}\rpa\; =\; \int_0^\infty e^{-ux} \lpa du - \sum_{n\ge 2} \delta_n (du) \rpa.$$
By (\ref{conv T}) and Theorem 7.1.1 in \cite{B}, this yields the following negative counterpart to Theorem \ref{GGCID}.

\begin{corollary}
\label{NoGG}
There exists $\a_0 < 1$ such that for every $\a\in (\a_0, 1),$ the law of $\X_\a$ is not a {\em GGC}.
\end{corollary}
This also implies that there is a function $\delta\colon (\a_0,1) \to [0,1)$ such that $\X_{\a,\rho}$ is not a GGC for $\a \in (\a_0,1)$ and $\rho \in [\delta(\a),1]$. 
Observe on the other hand that it does not seem possible to apply our methods to $\X_{\a,\rho}$ with a fixed $\rho\in (0,1).$ Indeed, as in the classical case, the possible limit laws of affine transformations of $\X_{\a,\rho}$ with $\rho\in (0,1)$ fixed and $\a\to 1$ are given only in terms of $\X_{1,\rho},$ whose law is a GGC.}
\end{remark}

\subsection{Proof of Theorem 4}

\subsubsection{The one-sided case} By (\ref{Kant1}) and Corollary 3.2 in \cite{JS}, we have the independent factorisation 
$$\X_\a\;\elaw\; b_\a\,\U^{-1/\beta}(1\, +\, \X)$$
where $1/\beta = 1/\a - 1$ and $\X$ has a CM density on $(0,\infty).$ We will now show the WS property for all positive random variables of the type 
$$\Y\; =\; \U^{-1/\beta}(1\, +\, \X)\, -\, 1$$
with $\beta > 0$ and $\X$ having a CM density on $(0,\infty).$ Setting $f,g$ for the respective densities of $\Y, \X,$ the multiplicative convolution formula shows that
$$f(x)\; =\; \frac{\beta}{(x+1)^{\beta +1}} \int_0^x (y+1)^\beta\, g(y)\, dy\; =\; \frac{\beta}{(1+1/x)^{\beta +1}} \int_0^1 (y+1/x)^\beta\, g(xy)\, dy$$
for every $x >0.$ In particular, one has $f(0+) = f(+\infty) = 0.$ Moreover, the first equality and an induction on $n$ imply that $f$ is smooth with
\begin{equation}
\label{Indu}
(x+1)f^{(n+1)}(x)\; =\; \beta g^{(n)}(x)\, -\, (\beta + n +1)f^{(n)}(x)
\end{equation}
for every $n\ge 0.$ Hence, we also have $f^{(n)} (+\infty) = 0$ for all $n \geq 0$ and a successive application of Rolle's theorem yields 
$$\# \{ x \in (0, \infty) \; \vert \; f^{(n)} (x) = 0  \}\; \geq\; 1$$
for every $n\ge 1.$ Fix now $n \geq 1$ and suppose that there exist $0 < x_n^{(1)} < x_n^{(2)} < \infty$ such that 
$$f^{(n)} (x^{(1)}_n) \;=\; f^{(n)} (x^{(2)}_n) \;=\;  0.$$
By (\ref{Indu}) and the complete monotonicity of $g,$ we have
$$(-1)^n f^{(n+1)} (x^{(i)}_n)\; >\; 0$$
for $i=1,2.$ An immediate analysis based on the intermediate value theorem shows then that there must exist $x_n^{(3)}\in (x_n^{(1)}, x_n^{(2)})$ with 
$$f^{(n)} (x^{(3)}_n) \;=\;  0\qquad\mbox{and}\qquad(-1)^n f^{(n+1)} (x^{(3)}_n)\; \le\; 0,$$
which is impossible again by (\ref{Indu}) and the complete monotonicity of $g.$ All in all, we have proved that 
$$ \# \{ x \in (0, \infty) \; \vert \; f^{(n)} (x) = 0  \}\; = \; 1$$ 
for all $n \geq 1,$ which is the WS property.

\qed

\subsubsection{The two-sided case} We know by Proposition A.1.4 in \cite{BerPat} that $f_{\a,\rho}$ is an analytic integrable function on $\rl$, and by Theorem 1.7 in \cite{HK} that it converges to zero at $\pm\infty,$ decreases near $+\infty$ and increases near $-\infty.$ Moreover, we have shown in Theorem 2 that if $\a\le 3/4,$ it is the density of an ID distribution on $\rl$ with L\'evy measure $\varphi_{\a,\rho}(x)\, dx$ such that $x\varphi_{\a,\rho}(x)$ and $x\varphi_{\a,\rho}(-x)$ are CM on $(0,\infty).$ We are hence in position to apply Corollary 1.2 in \cite{Kw}, which shows that $f_{\a,\rho}$ is BS.

\qed

\subsubsection{The exceptional 1-stable case} We use the second identity in Proposition \ref{prop:SF}, which rewrites
$$1\, -\,\F\;\elaw\; (1\,-\,\log\W)\; +\; \L.$$
We have seen in Remark \ref{MellWbis} (b) that the random variable $1-\log\W$ has a CM density on $(0,+\infty),$ in other words that it belongs to the class ${\rm ME}^*$ with the notations of \cite{TSPAMS}. Applying the Proposition in \cite{TSPAMS} with $n=1$ shows that $1-\F$ has a ${\rm WBS}_0$ density, with the notation of the main definition in \cite{TSPAMS}. As mentioned in the introduction, this means that the density of $\F$ is ${\rm WS}_-.$

\qed

\subsubsection{The two-sided 1-stable case with $b=0$ or $ab^{-1}\in\pi\ZZ$} We may suppose $a > 0$ by symmetry. If $b = 0$ the statement is clear since it is elementally shown that the Cauchy density
$$\frac{1}{\pi (1+x^2)}$$
is BS - see also Corollary 1.3 in \cite{Kw}. If $b\neq 0,$ we may suppose $b<0$ by symmetry. By independence, we have
$$\log \esp[e^{-\i\xi\C_{a,b}}]\; =\; -a\vert\xi\vert \; +\;\log\esp[e^{\i\vert b\vert\xi\F}],\qquad\xi\in\rl.$$
A further computation using Lemma 14.11 in \cite{Sato} and Remark \ref{MellWbis} (b) yields
$$\log \esp[e^{-\i\xi\C_{a,b}}]\; =\; c_1\; +\; c_2\i\xi\; +\;\int_\rl \lpa\frac{1}{\i\xi +s} -\lpa \frac{1}{s} \, -\, \frac{\i\xi}{s^2}\rpa\Un_{\rl\backslash (-1,1)} (s)\rpa \varphi_{a,b}(s)\, ds$$
for some $c_1,c_2\in\rl$ and
$$\varphi_{a,b} (s)\; =\; \frac{a}{\pi}\, s \,+ \,\lpa\vert b\vert s - [\vert b\vert s-1]_+\rpa\Un_{\{s\ge 0\}}.$$ 
This function satisfies (1.1) and (1.2) in \cite{Kw} and is such that $s\varphi_{a,b}(s)\ge 0.$ Moreover, for $ab^{-1}\in\pi\ZZ$ the function $\varphi_{a,b}(s) -k$ changes its sign only once for every $k\in\ZZ.$ Finally, we know from Propositions A.1.3 and A.2.1 in \cite{BerPat} that the density of $\C_{a,b}$ is smooth, converges to zero at $\pm\infty,$ decreases near $+\infty$ and increases near $-\infty.$ We can hence apply Theorem 1.1 in \cite{Kw} and conclude the proof.

\qed

\begin{remark}
\label{ExcF1}
{\em (a) If the random variable $1 - \log\W$ had a ${\rm PF}_\infty$ density as $\L$ does, then the ${\rm BS}$ character of $f_{1,1/2}$ and the additive total positivity arguments used in \cite{Kw, TSPAMS} would show that $\C_{a,b}$ has a BS density on $\rl$ for $a\neq 0.$ But $1- \log\W$ cannot have a ${\rm PF}_\infty$ density, since its law is not a GGC - see e.g. Example 3.2.2 in \cite{B}. 

\vspace{2mm}

(b) If $ab^{-1}\not\in\pi\ZZ,$ the function $\varphi_{a,b}(s) -k$ changes its sign at least three times for every negative integer $k$, so that we cannot use Theorem 1.1 in \cite{Kw}. It is not clear to the authors whether  the density of $\C_{a,b}$ is always BS for $a\neq 0,$ and the case $ab^{-1}\in\pi\ZZ$ might be more the exception than the rule.}
\end{remark}

\section{Further results}
\label{Further}

\subsection{Some properties of the function $\theta_\a$} In this paragraph we consider further aspects of the function 
\begin{equation}
\label{ThFa}
\theta_\a(r)\; =\;\frac{1}{\pi}\,\arg[F_\a(re^{-\i\pi})],
\end{equation}
whose non-decreasing character amounts to the GGC property for the law of $\K_{1-\a}.$ We first prove the following asymptotic result.

\begin{proposition}
\label{AsympFa}
For every $\a \in [1/5, 1)$, one has
$$\lim_{r \rightarrow +\infty}\Re (F_\a(re^{-\i\pi}))\; =\; \frac{1}{2\alpha}\cdot$$
For every $\a \in (0,1/5)$, one has 
$$\liminf_{r \rightarrow +\infty}\Re (F_\a(re^{-\i\pi}))\, =\, - \infty \qquad \text{and} \qquad \limsup_{r \rightarrow +\infty}\Re (F_\a(re^{-\i\pi}))\, =\, +\infty.$$
\end{proposition}

The second part of this proposition has an immediate corollary, which answers in the negative an open problem stated in \cite{JS} - see Conjecture 3.1 therein.

\begin{corollary}
\label{NegKa}
The function $\theta_\a$ is not monotonous on $(0,\infty)$ for $\a < 1/5.$ In particular, the law of $\K_\a$ is not a {\em GGC} for $\a > 4/5.$
\end{corollary}
 
\medskip

\noindent
{\em Proof of Proposition \ref{AsympFa}}. We have seen during the proof of Proposition \ref{ThInc} that $\Re (F_{1/2}(re^{-\i\pi}))\equiv 1$ and that 
$$\Re (F_\a(re^{-\i\pi}))\;\to\; \frac{1}{2\a}$$
as $r\to +\infty$ for all $\a\in (1/2,1).$ We next consider the case $\a\in [1/5,1/2)$ introducing, as above, the function
$$G_\a(r)\; =\; \Re(F_\a(r^{\frac{\a}{\a-1}}e^{-\i\pi})) \; =\; 1\;+\; \frac{1}{2\pi}\,\Im\lpa\int_0^\infty e^{-ru+ u^\a e^{ 2\i\pi\alpha}}\,\frac{du}{u}\rpa.$$
Setting $\theta = \frac{5}{6}(1-2\a)\in (0,1/2],$ we have $2\alpha+ \alpha \theta\in[1/2, 1)$ and by Cauchy's theorem, we can rewrite
$$G_{\alpha}(r)\; = \;1\, +\, \frac{\theta}{2}\, +\, \frac{1}{2\pi}\,\Im\lpa \int_0^\infty e^{-rue^{i\pi\theta}+ u^\alpha e^{i\pi(2\alpha+ \alpha \theta)}}\,\frac{du}{u}\rpa.$$
The latter converges to 
$$1\, +\, \frac{\theta}{2}\, +\, \frac{1}{2\pi}\,\Im\lpa \int_0^\infty e^{u^\alpha e^{i\pi(2\alpha+ \alpha \theta)}}\,\frac{du}{u}\rpa\; =\;1\, +\, \frac{\theta}{2}\, +\, \frac{1}{2\pi\a}\,\Im\lpa \int_0^\infty e^{-u e^{-i\pi(1-2\alpha- \alpha \theta)}}\,\frac{du}{u}\rpa $$
as $r\to 0.$ The evaluation of the oscillating integral on the right-hand side is given e.g. in Formula 1.6(36) p.13 in \cite{EMOT}, and we finally obtain
$$\lim_{r \rightarrow 0}G_{\alpha}(r)\; =\; 1 +  \frac{\theta}{2} + \frac{1}{2\alpha} (1- 2\a- \alpha\theta)\; =\; \frac{1}{2\alpha}\cdot$$

\medskip

We finally consider the case $\alpha \in (0,1/5)$, which is much more technical and requires several steps. Setting $\theta = 2\alpha/(1-\alpha) \in (0, 1/2)$, we have $2\alpha+ \alpha \theta =\theta $ and the same argument as above implies
$$G_{\alpha}(r)\; = \;1 \;+\; \frac{\theta}{2}\; +\; \frac{1}{2\pi}\,\Im\lpa \int_0^\infty e^{(-rt+t^\alpha)e^{\i\pi\theta}}\,\frac{dt}{t}\rpa.$$
Hence, we are reduced to show that
$$\liminf_{r\rightarrow 0} H_\alpha (r)\; = \; - \infty \qquad \text{and} \qquad \limsup_{x\rightarrow 0} H_\alpha (r)\; =\; + \infty$$
with the notations $f_r(t) = \sin(\pi\theta)(-rt +t^\a)$ and
$$H_\a(r)\; =\; \Im\lpa \int_0^\infty e^{(-rt+t^\alpha)e^{\i\pi\theta}}\,\frac{dt}{t}\rpa\; =\; \int_0^\infty e^{\cot (\pi\theta) f_r(t)}\,\sin(f_r(t))\,\frac{dt}{t}\cdot
$$
Let us begin with the liminf. Setting
$$r_k\; =\; \alpha \left(\frac{(1-\alpha)\sin(\pi\theta)}{2k\pi} \right)^{1/\alpha-1}\quad\mbox{and}\qquad m_k \;=\; \left( \frac{2k\pi}{(1-\alpha)\sin(\pi\theta)}\right)^{1/\alpha},$$
it is clear that the function $f_{r_k}\!(t)$ increases on $(0, m_k)$ and decreases on $(m_k, +\infty)$, and that its global maximum equals $f_{r_k}\!(m_k) = 2k\pi$. This yields 
$$\int_{m_k}^\infty e^{\cot (\pi\theta) f_{r_k}\!(t)}\,\sin(f_{r_k}\!(t))\,\frac{dt}{t}\; <\; 0\qquad\mbox{for every $k\ge 1.$}$$ 
Considering now the unique $a_k \in (0,m_k)$ such that $f_{r_k}(a_k) = \pi $, we have $\lim_{k\rightarrow \infty} a_k = (\pi / \sin(\pi\theta))^{1/\alpha}$, so that 
$$\int_0^{a_k} e^{\cot(\pi\theta) f_{r_k}\!(t)}\,\sin(f_{r_k}\!(t))\,\frac{dt}{t}\;\to\;\int_0^{ (\pi / \sin(\pi\theta))^{1/\alpha}}\!\!\!\! e^{\cos(\pi\theta) t^\a}\,\sin(t^\a\sin(\pi\theta))\,\frac{dt}{t}\; <\; \infty$$
as $k\to +\infty.$ Hence it suffices to show that $A_k \to - \infty$ as $k\to +\infty,$ with
$$A_k\; =\; \int_{a_k}^{m_k} e^{\cot(\pi\theta) f_{r_k}\!(t)}\,\sin(f_{r_k}\!(t))\,\frac{dt}{t}\; =\;\frac{1}{\sin(\pi\theta)}\int_{\pi}^{2k\pi} \frac{ e^{\cot(\pi\theta) u}}{-r_k \varphi_k(u) + \alpha (\varphi_k(u))^{\alpha}} \,\sin(u)\, du,$$
where the second equality comes from a change of variable, having set $\varphi_k(u)$ for the inverse function of $f_{r_k}$ on $[\pi, 2k\pi]$ and written
$$\varphi_k'(u)\; =\; \frac{1}{f'_{r_k}(\varphi_k(u))}\; = \frac{1}{\sin(\pi\theta)(-r_k + \alpha (\varphi_k(u))^{\alpha - 1})}\; >\; 0.$$ 
We next define $p_k(u) :=e^{-\cot(\pi\theta)u} (-r_k \varphi_k(u) + \alpha (\varphi_k(u))^{\alpha})$ and prove its strict unimodality on $[\pi, 2 k\pi],$ computing
$$p'_k(u) =e^{-\cot(\pi\theta) u} \frac{\varphi'_k(u)}{\varphi_k(u)} (-r_k t + \alpha^2 t^\a - \cos(\pi\theta)(-r_k t + \alpha t^\alpha)^2)$$
with $t =\varphi_k(u).$ The strict unimodality of $p_k(u)$ on $(\pi, 2k\pi)$ amounts to the fact that
$$q_k(t) \; =\; -r_k t + \alpha^2 t^\alpha - \cos(\pi\theta)(-r_k t + \alpha t^\alpha)^2$$
has at most one zero point on $[a_k, m_k]$. It is clear by construction that there exists $c_k \in (0, m_k)$ such that $g_k(t) = -r_k t + \alpha t^\alpha$ increases on $(0, c_k)$ and decreases on $(c_k, m_k)$, and for all $t \in (c_k, m_k)$ we have $q_k(t) = t g'_k(t) - \cos(\pi\theta)(-r_k t + \alpha t^\alpha)^2 < 0.$ On the other hand, the function $g_k(t)$ is increasing and concave on $[0, c_k)$, so that its inverse function $\psi_k(v)$ is increasing and convex on $[0, g_k(c_k))$. Now since 
$$q_k(t)\, =\, 0\;  \Leftrightarrow\; \cos(\pi\theta) v^2 \,-\, \alpha v\, +\, (1-\alpha) r_k \psi_k(v)\, =\, 0,$$
we see that there are at most two solutions of $q_k(t) = 0$ on $[0,c_k),$ one of them being zero, and hence at most one solution on $[a_k,m_k),$ as required. We now denote by $z_k$ the unique mode of $p_k(u)$ on $[a_k,m_k]$ and, setting $l_k =\inf\{l\ge 1, \, z_k\le 2l\pi\},$ decompose
$$A_k \; =\; \frac{1}{\sin(\pi\theta)}\lpa \int_\pi^{2 l_k\pi}p_k^{-1}(u)\,\sin(u)\, du \; +\;\int_{2 l_k \pi}^{2 k\pi}p_k^{-1}(u)\,\sin(u)\, du \rpa.$$
Since $z_k\to\tan(\pi\theta)$ viz. $l_k\to l_\infty < +\infty$ as $k\to \infty,$ it is easy to see that the first term in the decomposition is bounded, and we are finally reduced to show that
$$B_k\; =\; \int_{2 l_k \pi}^{2 k\pi}p_k^{-1}(u)\,\sin(u)\, du\; \to\; -\infty\qquad\mbox{as $k\to +\infty.$}$$
Since $p_k^{-1}(u)$ increases on $[2 l_k\pi, 2k\pi],$ we have 
$$B_k\; =\; \sum_{j= l_k}^{k-1}\lpa\int_{2j\pi}^{(2j +1)\pi} (p_k^{-1}(u) -p_k^{-1}(u+\pi))\,\sin(u)\, du\rpa$$
for every $k\ge 1$ and since $p_k (u) \to \frac{\alpha}{\sin(\pi\theta)}u e^{-\cot(\pi\theta)u}$ pointwise as $k\to +\infty,$ Fatou's lemma implies
$$\limsup_{k\rightarrow +\infty} B_k \; \le \; \frac{\sin(\pi\theta)}{\alpha} \sum\limits_{j=l_\infty}^{\infty} \int_{2j\pi}^{(2j+1)\pi}\lpa\frac{e^{\cot(\pi\theta) u}}{u} - \frac{e^{\cot(\pi\theta) (u+\pi)}}{u+\pi} \rpa \sin (u)\, du.$$
Using the inequality 
\begin{eqnarray*}
\frac{1+ e^{\pi\cot\pi\theta}}{2u} \leq \frac{e^{\pi\cot\pi\theta}}{u+\pi} 
\end{eqnarray*}
which holds for $u \geq \frac{\pi(e^{\pi\cot\pi\theta}+1)}{e^{\pi\cot\pi\theta}-1}$,  we deduce that for $j_\infty$ large enough, one has
\begin{eqnarray*}
\limsup_{k\rightarrow +\infty} B_k & \le & \frac{\pi \sin(\pi\theta)}{\alpha} \sum\limits_{j=j_\infty}^{\infty} \frac{1-e^{\pi\cot\pi\theta}}{2}\int_{2j\pi}^{(2j+1)\pi}\frac{e^{\cot(\pi\theta)u} }{u} \sin(u)\, du \\
&\le & -\frac{(e^{\pi\cot\pi\theta}-1)\pi\sin(\pi\theta)}{4\alpha} \sum\limits_{j=j_\infty}^{\infty}  \int_{2j\pi + \pi /6 }^{2j\pi +5\pi /6}\frac{1}{u} \,du \; =\; -\infty. 
\end{eqnarray*}
All of this shows that 
$$\liminf_{r\rightarrow 0} H_\a(r) \;=\; - \infty.$$
The argument for the limsup follows exactly along the same lines, considering the subsequence
$$ \tilde{r}_k\; =\; \alpha \left(\frac{(1-\alpha)\sin(\pi\theta)}{(2k+1)\pi} \right)^{1/\alpha -1}.$$

\qed

\begin{remark}
\label{Asympbis}
{\em (a) In the case $\a\in [1/3, 1/2)$ we have seen in the proof of Proposition \ref{ThInc} that
$$\Re(F_\a(re^{-\i\pi}))\; = \; 1\; +\; \frac{\rho}{2}\,\pb[\Z_{\a,\rho}^{-\a} \leq r^{1-\a}]$$
with $\rho = 1/\a- 2,$ which does converge to $1/(2\a)$ as $r\to +\infty.$ In the case $\a\in[1/4,1/3),$ the proof of Proposition \ref{ThInc} shows that
$$\lim_{r\to 0} G_\a(r)\; =\; 1\; +\; \frac{1}{2\pi\a} \int_0^\infty e^{\cos(2\pi\a) u} \sin(\sin(2\pi\a)u)\, \frac{du}{u}\; =\; \frac{1}{2\a}$$
again by Formula 1.6(36) in \cite{EMOT}. The above contour argument is hence only necessary for $\a\in[1/5, 1/4).$ 

\vspace{2mm} 

(b) As mentioned in Remark \ref{Rho1} (a), the above proof shows that $x\mapsto\pb[\Z_\a^{-\frac{\alpha}{1-\alpha}} \geq x]$ is not HCM for every $\a < 1/5.$ By Theorem 6.3.5 in \cite{B}, this implies that $\Z_\a^{-\frac{\alpha}{1-\alpha}}$ is not HCM for $\a < 1/5$ either. This shows that Conjecture 1.2 in \cite{PB} is not true in general.} 
\end{remark}

We believe that $\theta_a$ is non-decreasing for $\a\in [1/5, 1),$ which is equivalent to the following

\begin{conjecture}
\label{KaGGC}
The law of $\K_\a$ is a {\em GGC} if and only if $\a\le 4/5.$
\end{conjecture}

\begin{figure}
\begin{center}
%\advance\rightskip-7cm
\begin{minipage}{0.45\hsize}
\begin{center}
%\advance\leftskip-5cm
\includegraphics[scale=0.2, bb=0 0 960 731]{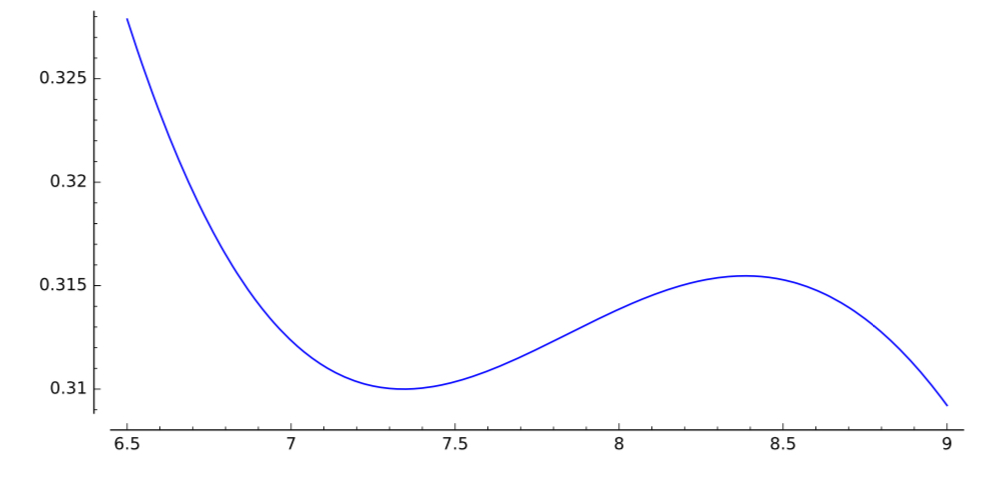} 
\caption{$x^{-1} G_{1/5}(x^{-5})$} \label{G0.2}
\end{center}
  \end{minipage}
\begin{minipage}{0.45\hsize}
\begin{center}
\includegraphics[scale=0.2, bb=0 0 960 731]{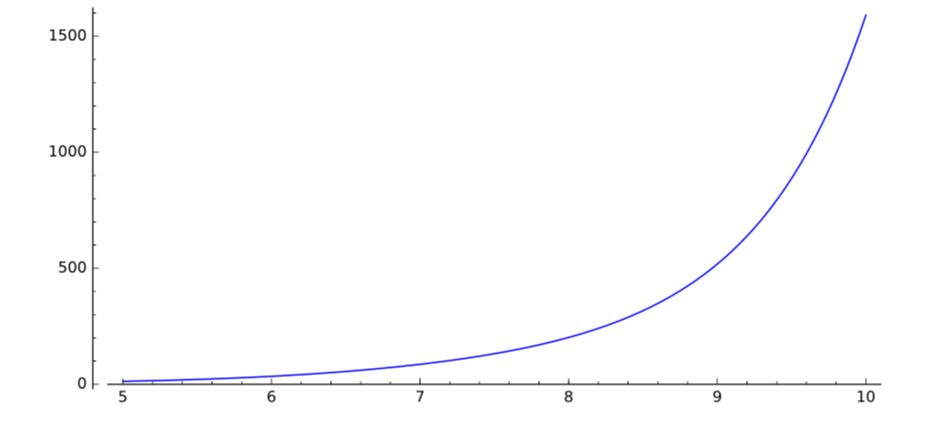} 
\caption{${\tilde G_{1/5}}(x)$} \label{Gt0.2}
\end{center}
\end{minipage}
\end{center}
\end{figure}

The above Corollary \ref{NegKa} shows the only if part, and in the proof of Theorem 2 we have shown the if part for $\a\le 3/4.$ However, it seems that our methods fail to handle the remaining case $\a\in(3/4,4/5],$ because some simulations show that $r^{\a-1}G_\a(r^{\frac{\a-1}{\a}})= r^{\a -1}\Re(F_\a(re^{-\i\pi}))$ is not monotonous anymore, at least for $\a$ close enough to 1/5 - see Figure \ref{G0.2}. Observe from (\ref{Wrgt}) that the problem can be reformulated in terms of the monotonicity of the ratio of two power series, the non-decreasing character of $\theta_\a$ being equivalent to that of
$${\tilde G}_\a \, :\, x\;\longmapsto\; \frac{\Im(F_\a(x e^{-\i\pi}))}{\Re(F_\a(x e^{-\i\pi}))}\;=\;\frac{\displaystyle \sum\limits_{n \geq 0} \frac{\sin (n\pi \alpha)}{n! \Gamma(1-n\alpha)}\,x^n}{\displaystyle \sum\limits_{n \geq 0} \frac{\cos (n\pi \alpha)}{n! \Gamma(1-n\alpha)}\, x^n}$$
on $(0,\infty).$ A necessary condition for ${\tilde G}_\a$ to be non-decreasing is that its denominator does not vanish on $(0,\infty)$, which is false for $\a < 1/5$ by Proposition \ref{AsympFa} and true for $\a\ge 1/4$ by the proof of Theorem 2. But the case $\a\in[1/5,1/4)$ still eludes us. Let us mention that monotonicity properties of ratios of power series are studied in the literature on special functions - see e.g. Chapter 3.1 in \cite{AB}. For example, one could be tempted to apply Theorem 4.3 in \cite{HVV} since $x\mapsto\tan(x\pi\a)$ is locally increasing. However, we could not find any clue in this literature for our problem, and it is not easy to understand why the value $\a =1/5$ should be critical for the monotonicity of the above ratio. See Figure 4 for a convincing simulation. Let us finally mention \cite{Ne} for an operator-theoretic approach to the above power series.\\

We finally turn to the behaviour of $F_\alpha(re^{-\i\pi})$ at infinity, which implies that of $\theta_\a(r).$ 

\begin{proposition}
\label{Asymp}
One has
$$F_\alpha(re^{-\i\pi}) \;\sim\; \frac{\i\,c_\a e^{b_{1-\a} r}}{\sqrt{r}} \qquad\mbox{as $r \rightarrow +\infty,$}$$
with $c_\a =\frac{\a^{\frac{1}{2(\a -1)}}}{\sqrt{2\pi(1-\a)}}.$ In particular, one has $\theta_\a(r)\to 1/2$ as $r\to +\infty.$
\end{proposition}

\proof

From (\ref{Wrgt}), we can write
$$F_\alpha(re^{-\i\pi})\; =\; \phi (-\a, 1, r^{1-\a} e^{\i\pi\a}),\qquad r >0.$$
We now use the asymptotic expansion for large $z\in\CC$ and $a\in (-1,0)$ of the Wright function $\phi(a,b,z),$ which has been obtained in \cite{W}. Applying therein Theorem 1 for $\a\le 1/3$ resp. Theorem 5 for $\a > 1/3$ and taking the first term in (1.3) implies the required asymptotic for $F_\a(re^{-\i\pi}),$ since we have here
$$A_0 \; =\; \frac{1}{\sqrt{2\pi\a}}\qquad\mbox{and}\qquad Y \; =\; b_{1-\a} re^{-\i\pi}$$
in the notation of \cite{W}, the first equality being a consequence of Stirling's formula. From (\ref{ThFa}), we then readily deduce that $\theta_\a(r)\to 1/2$ as $r\to +\infty.$

\endproof

\begin{remark}
\label{Asympbi}
{\em (a) Taking the first two terms in the series representation (\ref{Wrgt}) yields at once the asymptotic behaviour of $\theta_\a(r)$ at zero, which is
$$\theta_\a(r)\;\sim\; \frac{r^{1-\a}}{\Ga(\a)\Ga(1-\a)^2}\cdot$$
On the other hand, the complete asymptotic expansion (1.3) in \cite{W} has only purely imaginary terms in our framework, so that we cannot deduce from it the asymptotics of $\theta_\a(r) -1/2$ at infinity. It follows from Proposition \ref{ThInc} that $\theta_\a(r)\in [0,1/2)$ for $\a\ge 1/4,$ and from Proposition \ref{AsympFa} that $\theta_\a(r) -1/2$ crosses zero an infinite number of times for $\a <1/5,$ as $r\mapsto +\infty.$ For $\a\in[1/5, 1/4),$ we are currently unable to prove that $\theta_\a(r)\in [0,1/2)$ for every $r >0,$ which would be a first step to show that it increases from 0 to 1/2. Recall that the latter is equivalent to the fact that the denominator of the above ${\tilde G}_\a$ does not vanish on $(0,\infty).$

\vspace{2mm}

(b) If $\a\le 3/4,$ it follows from (\ref{RepFa}), Theorem 8.2 and Remark 8.3 in \cite{SSV}, and the above proposition, that the Thorin mass of the GGC random variable $\K_\a$ equals $1/2.$ Hence, $\K_\a-b_\a$ is a $\Ga_{1/2}-$mixture by Theorem 4.1.1. in \cite{B}, which is a refinement of Corollary 3.2 in \cite{JS}. Since this property amounts to the CM character of $x\mapsto \sqrt{x}\,f_{\K_\a}(b_\a +x),$ a perusal of the proof of Theorem \ref{freeID} shows that $\X_\a - b_\a$ is a $\Ga_{3/2}-$mixture as soon as $\a\le 3/4.$ We believe that this is true for every $\a\in (0,1).$ }
\end{remark}

\subsection{An Airy-type function} In this paragraph, we discuss a curious connection between the two cases $\a =1/3$ and $\a =1/5$ in the analysis of the function 
$$G_\a(r) \; =\; 1\;+\; \frac{1}{2\pi}\,\Im\lpa\int_0^\infty e^{-ru+ u^\a e^{ 2\i\pi\alpha}}\,\frac{du}{u}\rpa.$$
The latter was important during the proofs of Theorem 2 and Proposition \ref{AsympFa}. For $\a =1/3,$ a contour integration as in Proposition \ref{AsympFa} with $\theta = -1/2$ implies, making the change of variable $s = (3r)^{-1/3},$
\begin{eqnarray*}
G_{\frac{1}{3}} (r) & = & \frac{3}{4} \; +\;  \frac{3}{2\pi}\,\int_0^\infty \,\sin (u^3/3 + us)\, \frac{du}{u}\\
& = & \frac{3}{4} \; +\;  \frac{3}{2\pi}\lpa\int_0^\infty \sin(u^3/3)\,\frac{du}{u}\; +\; \int_0^\infty \lpa\int_0^s\,\cos (u^3/3 + u z) \, dz\rpa du\rpa\\
& = & 1 \; +\;  \frac{3}{2\pi}\int_0^s \lpa\int_0^\infty\,\cos (u^3/3 + u z) \, du\rpa dz\\
& = & 1 \; +\;  \frac{3}{2}\int_0^s {\rm Ai}(z)\, dz
\end{eqnarray*}
where ${\rm Ai}$ stands for the classic Airy function - see e.g. Paragraph 7.3.7 in \cite{EMOT}. In particular, we retrieve the fact that
$$r\;\mapsto\;(3r)^{\frac{1}{3}}G_{\frac{1}{3}} (r)\; =\; \frac{1}{s}\; +\;  \frac{3}{2}\int_0^1 {\rm Ai}(sz)\, dz$$
increases, by the well-known decreasing character of Ai on $(0,\infty).$ For $\a = 1/5,$ the contour integration of Proposition \ref{AsympFa} with $\theta = 1/2$ yields, with the change of variable $s = (5r)^{-1/5},$
\begin{eqnarray*}
G_{\frac{1}{5}} (r) & = & 1 \; +\;  \frac{5}{2}\int_0^s \lpa\int_0^\infty\,\cos (u^5/5 - u z) \, du\rpa dz\\
& = & 1 \; +\;  \frac{5}{2}\int_0^s {\rm Ai}^{(5)}(-z)\, dz
\end{eqnarray*}
where we have defined, for every integer $k\ge 3,$ the semi-converging integral
$${\rm Ai}^{(k)}(x)\; =\;\frac{1}{\pi}\int_0^\infty\,\cos (u^k/k + u x) \, du, \qquad x\in\rl.$$
We did not find any reference on the above Airy-type functions in the literature, which are solution to some linear ODE of higher order. Observe that similarly as above, one has
$$(5r)^{\frac{1}{5}}G_{\frac{1}{5}} (r)\; =\; \frac{1}{s}\; +\;  \frac{5}{2}\int_0^1 {\rm Ai}^{(5)}(-sz)\, dz$$
but here we cannot deduce any conclusion on the monotonicity of $r^{\frac{1}{5}}G_{\frac{1}{5}} (r)$ because of the negative sign in the Airy-type function. The simulation displayed in Figure 5 shows indeed that ${\rm Ai}^{(5)}(-x)$ exhibits on $(0,\infty)$ exactly the same damped oscillating behaviour as ${\rm Ai}(-x).$ It could be interesting for our purposes to perform a rigorous study of the functions ${\rm Ai}^{(k)},$ as in the case $k=3$ with the Bessel functions. We leave this analysis for future research.

\begin{figure}
\begin{center}
\includegraphics[scale=0.3, bb=0 0 760 470]{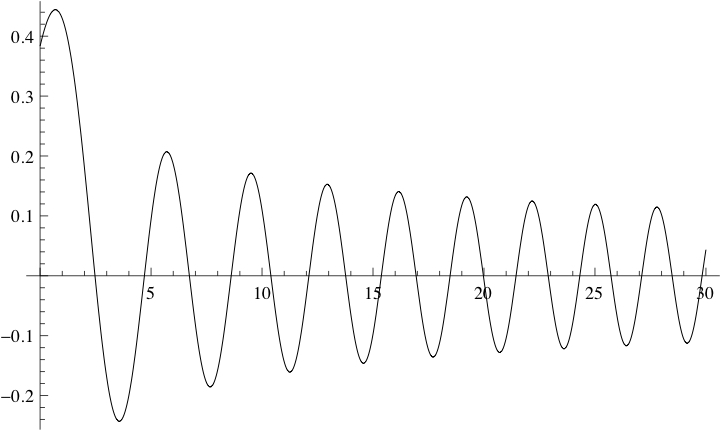}
\caption{${\rm Ai}^{(5)}(-x)$} \label{Ai5}
\end{center}
\end{figure}

\subsection{Asymptotic expansions for the free extreme stable densities} In this paragraph we derive the full asymptotic expansion at zero of the density $f^{}_{\Y_\a}$ of the random variable
$$\Y_\a\; =\; \lacc\begin{array}{ll}
 \X_\a - b_\a & \mbox{if $\a\in (0,1),$}\\
\X_{\a,1-1/\a} + b_{1/\a}^{-1/\a} \,\elaw\, b_{1/\a}^{-1/\a} - \X_{\a,1/\a} & \mbox{if $\a\in (1,2],$}
\end{array}\right.$$
and $\Y_1 = 1-\F.$ We will use the standard notation of Definition C.1.1 in \cite{AAR} for asymptotic expansions. Our expansions complete the estimates of Proposition A.1.2 in \cite{BerPat} and the series representations of Theorem 1.7 in \cite{HK}, from which one can only infer that the random variable $\Y_\a$ is positive. They can also be viewed as free analogues of Linnik's expansions (14.35) in \cite{Sato} - see also Theorem 2.5.3 in \cite{Z} - for the classical extreme stable distributions. Observe that in the classical case, the expansion for $\a >1$ is deduced from that of the case $\a\in[1/2,1)$ by the Zolotarev's duality which is discussed in Section 2.3 of \cite{Z}. Even though the very same duality relationship holds in the free case - see Proposition A.3.1 in \cite{BerPat} and Corollary 1.4 in \cite{HK}, for $\Y_\a$ this duality only yields
$$f^{}_{\Y_{1/\a}}(x)\; =\; \frac{1}{\a}\,(b_\a^{-\a} -x)^{-1/\a-1}f^{}_{\Y_\a}((b_\a^{-\a} -x)^{1/\a} -b_\a)$$
for every $\a\in[1/2,1),$ and does not seem particularly helpful to connect explicitly the two expansions at zero. 
When $\a\neq 1,$ our method hinges on Wright's original papers \cite{W0} for the case $\a > 1$ and \cite{W} for the case $\a < 1.$ It is remarkable that the two expansions turn out to have the same parametrization.

\begin{proposition}
\label{W0a}
For every $\a\in(0,1)\cup(1,2],$ one has 
$$f_{\Y_\a}(x)^{}\; \sim\; \sum_{n=0}^\infty a_n(\a)\, x^{n+1/2}\qquad\mbox{as $x\to 0,$}$$ 
with 
$$a_n(\a)\; =\; \lpa\frac{2}{\a}\rpa^{n+1/2}\!\!\!\frac{(-1)^n}{\pi\,\vert\a-1\vert^{(n+3/2)/\a}\,(2n+1)!}\,\times\,\frac{{\rm d}^{2n}}{{\rm d} v^{2n}}\lpa (1-v)^{-2}\pFq{2}{1}{\a+1,1}{3}{v}^{-n-1/2}\rpa_{v=0}.$$ 
\end{proposition}

\proof We begin with the case $\a > 1,$ writing down first $f_{\Y_\a}$ with the help of Bromwich's integral formula
$$f_{\Y_\a}^{}(x)\; =\;\frac{1}{2\pi\i} \int_{1-\i\infty}^{1+\i\infty} e^{zx}\cL_\a(z)\, dz,$$
where
$$\cL_\a(z)\; =\; \esp[e^{-z\Y_\a}]\; =\; e^{-\a(\a-1)^{\frac{1}{\a} -1}z}\,\times\,\esp[e^{z \X_{\a, 1/\a}}]$$
is well-defined and analytic on the open right half-plane. Combining next Theorem 1.8 in \cite{HK} and Theorem 2 in \cite{W0}, we obtain
$$\cL_\a(z)\; =\; e^{-\a(\a-1)^{\frac{1}{\a} -1}\,z}\,\times\,\phi(\a-1,2,z^\a)\; =\; {\rm O} (\vert z\vert^{-3/2})$$
uniformly on the right half-plane. Making a change of variable and applying Cauchy's theorem, we deduce
$$f_{\Y_\a}^{}(x)\; =\;\frac{1}{2\pi\i\, x} \int_{x-\i\infty}^{x+\i\infty} e^{z}\cL_\a(zx^{-1})\, dz\; =\;\frac{1}{2\pi\i\, x} \int_{1-\i\infty}^{1+\i\infty} e^{z}\cL_\a(zx^{-1})\, dz.$$
Using now the full asymptotic expansion of Theorem 2 in \cite{W0}, we get
$$f_{\Y_\a}(x)^{}\; \sim\; \sum_{n=0}^\infty a_n(\a)\, x^{n+1/2}\qquad\mbox{as $x\to 0,$}$$ 
where 
\begin{eqnarray*}
a_n(\a) &=& \frac{(-1)^n\,a_n}{(\a-1)^{(n+3/2)/\a}}\lpa\frac{1}{2\pi\i} \int_{1-\i\infty}^{1+\i\infty} e^{z} z^{-3/2-n}\, dz\rpa\; = \;  \frac{(-1)^n\,a_n}{(\a-1)^{(n+3/2)/\a}\Ga(n+3/2)}
\end{eqnarray*}
and $a_n$ is defined at the beginning of p.258 in \cite{W0} for $\rho =\a-1$ and $\beta = 2.$ Above, the interchanging of the contour integral and the expansion is easily justified - alternatively one can use the generalized Watson's lemma which is mentioned at the top of p.615 in \cite{AAR}, whereas the second equality follows from Hankel's formula - see e.g. Exercise 1.22 in \cite{AAR}. To conclude the proof of the case $\a > 1,$ it remains to evaluate the coefficients $a_n(\a),$ which is done in observing that the function in (1.21) of \cite{W0} is here
$$\sqrt{\pFq{2}{1}{\a+1,1}{3}{v}},$$
and making some simplifications.

We now consider the case $\a < 1.$ The argument is analogous but it depends on the expansions of \cite{W} which, the author says, cannot be simply deduced from those of \cite{W0}. We again write
$$f_{\Y_\a}^{}(x)\; =\;\frac{1}{2\pi\i} \int_{1-\i\infty}^{1+\i\infty} e^{zx}\cL_\a(z)\, dz,$$
where
$$\cL_\a(z)\; =\; e^{\a(1-\a)^{\frac{1}{\a} -1}z}\,\times\,\esp[e^{-z \X_\a}]\; =\;e^{\a(1-\a)^{\frac{1}{\a} -1}\,z}\,\times\,\phi(\a-1,2,-z^\a)\; =\; {\rm O} (\vert z\vert^{-3/2})$$
uniformly in the open right half-plane, the second equality following from Theorem 1.8 in \cite{HK} and the estimate from the Lemma p.39 in \cite{W}. Reasoning as above, we get
$$f_{\Y_\a}(x)^{}\; =\;\frac{1}{2\pi\i\, x} \int_{1-\i\infty}^{1+\i\infty} e^{z}\cL_\a(zx^{-1})\, dz\; \sim\; \sum_{n=0}^\infty a_n(\a)\, x^{n+1/2}\qquad\mbox{as $x\to 0,$}$$ 
where 
\begin{eqnarray*}
a_n(\a) &=&  \frac{a_n}{(\a(1-\a)^{\frac{1}{\a} -1})^{n+3/2}\,\Ga(n+3/2)}
\end{eqnarray*}
and $a_n$ is defined at the bottom of p.38 in \cite{W} for $\sigma =1-\a$ and $\beta = 2.$ After some simplifications, we also obtain the required expression for $a_n(\a).$

\endproof
 
\begin{remark}
\label{W0ab}
{\em (a) It does not seem that a simple closed formula can be obtained for the coefficients $a_n(\a)$ in general. We can compute
$$a_0(\a)\; =\; \sqrt{\frac{2}{\a}}\,\times\, \frac{1}{\pi \vert\a-1\vert^{3/(2\a)}}\quad\mbox{and}\quad a_1(\a)\; = \; -\sqrt{\frac{2}{\a}}\,\times\,\lpa\frac{2\a^2 - 23\a +47}{36\pi\,\a\,\vert\a-1\vert^{5/(2\a)}}\rpa.$$
Observe that $a_1(\a)$ is always negative. We believe that in general, one has
$$a_n(\a)\; =\;\sqrt{\frac{2}{\a}}\,\times\, \frac{Q_{2n}(\a)}{\pi\, \a^n\,\vert\a-1\vert^{(2n+3)/(2\a)}}$$
for some $Q_{2n}\in\QQ_{2n}[X].$ This would again mimic the classical situation, save for the fact that here the polynomial $Q_{2n}$ does not seem to have symmetric coefficients - see Remark 2 p.101 in \cite{Z}.\\

(b) For $\a =2,$ the involved hypergeometric function becomes the standard geometric series and we simply get
$$a_n(2)\; =\; \frac{(-1)^n}{\pi\,(2n+1)!}\,\times\,\frac{{\rm d}^{2n}}{{\rm d} v^{2n}}\lpa (1-v)^{n-3/2}\rpa_{v=0}\; =\; \frac{-1}{\pi\,(2n-1)\, 16^n}\, \binom{2n}{n},$$
which is always negative except for $n=0.$ Of course, this can be retrieved via the binomial theorem for the explicit density
$$f_{\Y_2}(x)^{}\; =\; \frac{\sqrt{x}}{\pi}\,\sqrt{ 1 - \frac{x}{4}}\cdot$$

(c) For $\a =1/2,$ the involved hypergeometric function simplifies with the help of Exercise 3.39 in \cite{AAR}, and we get
$$a_n(1/2)\; =\; \frac{(-1)^n 4^{n+2}}{2\pi\,(2n+1)!}\,\times\,\frac{{\rm d}^{2n}}{{\rm d} v^{2n}}\lpa (1+\sqrt{1-v})^{2n+1}(1-v)^{-2}\rpa_{v=0}\; =\; \frac{(-1)^n (n+1)\, 4^{n+2}}{\pi},$$
whose signs alternate. This again can be retrieved via the binomial theorem for the explicit density
$$f_{\Y_{1/2}}(x)^{}\; =\; \frac{16\,\sqrt{x}}{\pi(1+4x)^2}\cdot$$

(d) As already observed in Remark \ref{KantFuss} (a), the densities of $\Y_{1/3}$ and $\Y_{2/3}$ can be written in closed form with the help of formul\ae\, (40) and (41) in \cite{MPZ}. In principle, a full asymptotic expansion can also be derived from these expressions, but the task seems too painful. Notice that here, the involved hypergeometric functions do not seem to simplify.\\

(e) The above proof shows that the following functions
$$\lbd\mapsto e^{-\a(\a-1)^{\frac{1}{\a} -1}\,\lbd}\,\phi(\a-1,2,\lbd^\a)\qquad\mbox{resp.}\qquad \lbd\mapsto e^{\a(1-\a)^{\frac{1}{\a} -1}\,\lbd}\,\phi(\a-1,2,-\lbd^\a)$$ 
on $(0,\infty)$, which are obtained in removing Wright's exponential term at infinity, are CM functions for $\a\in(1,2]$ resp. for $\a\in (0,1).$ \\

(f) For $\a > 1,$ we can also compute the Mellin transform of $\Y_\a,$ starting from the formula
$$\esp[\Y_\a^{-s}]\; =\; \frac{1}{\Ga(s)} \int_0^\infty \esp [e^{-\lbd \Y_\a}]\, \lbd^{s-1}\, d\lbd$$
which is valid for every $s > 0$ with possible infinite terms on both sides. This becomes here
$$\esp[\Y_\a^{-s}]\; =\; \frac{b_\a^{-s}}{\Ga(s)}\,\sum_{n\ge 0} \frac{\Ga(s+\a n)}{n!\, \Ga(2 + (\a-1) n)}\; b_\a^{-\a n}$$
with the notation $b_\a = \a(\a-1)^{\frac{1}{\a} -1}$ and has, by Stirling's formula, an analytic extension for $-\a < s < 3/2.$ Formally, this rewrites
$$\esp[\Y_\a^{s}]\; =\; \frac{b_\a^{-s}}{\Ga(s)}\;_1\Psi_1\lcr\left.\begin{array}{l} (-s,\a)\\ (2,\a-1)\end{array} \,\rva\, b_\a^{-\a}\rcr, \qquad -3/2 <s <\a,$$
where $_r\Psi_s$ is the generalized hypergeometric function originally studied in \cite{Fox,W1}, which is sometimes coined as a generalized Wright function, and which should not be confused with the $_r\psi_s$ hypergeometric series defined in (10.9.4) of \cite{AAR}. For $\a =2,$ Gauss's multiplication and summation formul\ae\, for the Gamma and the hypergeometric function - see Theorems 1.5.1 and 2.2.2 in \cite{AAR}, respectively - transform this expression into
$$\esp[\Y_2^{-s}]\; =\;2^s\, \pFq{2}{1}{-s/2,(1-s)/2}{2}{1}\; =\; \frac{4^{s+1}}{\sqrt{\pi}}\,\times\,\frac{\Ga(3/2 +s)}{\Ga(3+s)},$$
in accordance with $\Y_2\elaw 4 \B_{3/2,3/2}.$ }
\end{remark}  
We now complete the picture and derive the asymptotic expansion of $\Y_1 = 1 -\F.$ To state our result, we need to introduce the Stirling series $\{c_n,\, n\ge 0\}$ appearing in the expansion
$$\lpa \frac{\e}{x}\rpa^x \sqrt{\frac{x}{2\pi}}\;\Gamma(x)\;\sim\; \sum_{n\ge 0} c_n x^{-n}\qquad \mbox{as $x\to +\infty,$}$$
which is given e.g. in Exercise 23 p.267 of \cite{Co} - see also Lemma 1 in \cite{Fox} . One has $c_0 = 1, c_1 = 1/12, c_2 = 1/288$ and $c_3 = -139/51840.$ In general, $c_n$ is a rational number and the corresponding sequences of numerators and denominators are referenced under A00163 and A00164 in the online version of \cite{OEI}.

\begin{proposition}
\label{W1a}
One has
$$f_{\Y_1}(x)^{}\; \sim\; \sum_{n=0}^\infty a_n(1)\, x^{n+1/2}\qquad\mbox{as $x\to 0,$}$$ 
with 
$$a_n(1)\; =\;  \frac{(-1)^n 2^{2n+1/2}\, n!\,(c_0+\cdots + c_n)}{\pi\,(2n+1)!}\cdot$$ 
\end{proposition}

\proof Applying Remark \ref{MellWbis} (a), we first compute the Laplace transform
$$\esp[e^{-z\Y_1}]\; =\; \frac{1}{z(1+z)}\,\lpa \frac{z}{\e}\rpa^z\,\frac{1}{\Ga(z)}$$
for every $z$ in the open right half-plane. Comparing next (2.15) and (2.21) in \cite{NW}, we get the expansion
\begin{eqnarray*}
\esp[e^{-z\Y_1}] & \sim & \frac{1}{\sqrt{2\pi z} (1+z)}\,\sum_{n\ge 0} (-1)^n c_n\, z^{-n}\\
& \sim &  \frac{1}{\sqrt{2\pi z^3}}\,\lpa\sum_{n\ge 0} (-1)^n z^{-n}\rpa\lpa\sum_{n\ge 0} (-1)^n c_n\, z^{-n}\rpa\\
& \sim & \frac{1}{\sqrt{2\pi z^3}}\,\sum_{n\ge 0} (-1)^n (c_0+\cdots +c_n)\,z^{-n}.
\end{eqnarray*}
uniformly in the open right half-plane. Reasoning as in Proposition \ref{W0a}, we finally obtain
$$f_{\Y_1}(x)^{}\; \sim\; \sum_{n=0}^\infty a_n(1)\, x^{n+1/2}\qquad\mbox{as $x\to 0,$}$$ 
with 
$$a_n(1)\; =\; \frac{(-1)^n (c_0+\cdots +c_n)}{\sqrt{2\pi}\,\Ga(n+3/2)}\; =\; \frac{(-1)^n 2^{2n+1/2}\, n!\,(c_0+\cdots + c_n)}{\pi\,(2n+1)!}\cdot$$ 
\endproof

\begin{remark}
\label{W1ab}
{\em It is easy to see from (\ref{conv T}) that
$$(1-\a)^{-1}\Y_\a\; \claw\; \Y_1\qquad\mbox{as $\a\uparrow 1,$}$$
and it is natural to infer from this and Proposition \ref{W0a} that
\begin{eqnarray*}
a_n(1) & = & \frac{(-1)^n \, 2^{n+1/2}}{\pi\,(2n+1)!}\,\times\,\frac{{\rm d}^{2n}}{{\rm d} v^{2n}}\lpa (1-v)^{-2}\pFq{2}{1}{2,1}{3}{v}^{-n-1/2}\rpa_{v=0}\\
& = & \frac{(-1)^n}{\pi\,(2n+1)!}\,\times\,\frac{{\rm d}^{2n}}{{\rm d} v^{2n}}\lpa \frac{v^{2n+1}}{(1-v)^2(-v-\log(1-v))^{n+1/2}}\rpa_{v=0},
\end{eqnarray*} 
except that we cannot interchange a priori the asymptotic expansion at zero and the convergence in law. We have checked the correspondence for $n = 0$ and $n=1$, with
$$a_0(1)\; =\; \frac{\sqrt{2}}{\pi}\qquad\mbox{and}\qquad a_1(1)\; =\; -\frac{13\sqrt{2}}{18\pi}$$
to be compared with Remark \ref{W0ab} (a). We believe that this formula is true for every $n\ge 1.$ Observe that this is equivalent to the following expression of the Stirling series:
$$c_n\; =\; b_n\; -\; b_{n-1}, \qquad n\ge 1,$$
with 
$$b_n\; =\; \frac{1}{2^{2n+1/2}\,n!}\,\times\,\frac{{\rm d}^{2n}}{{\rm d} v^{2n}}\lpa \frac{v^{2n+1}}{(1-v)^2(-v-\log(1-v))^{n+1/2}}\rpa_{v=0},$$

\smallskip

\noindent
which is different from the combinatorial expression given in Exercise 23 p.267 of \cite{Co}, and which we could not locate in the literature.}
\end{remark}
\subsection{Product representations for  $\K_\a$ and $\X_\a$ with $\a$ rational} In the classical framework, the following independent factorization of the positive stable random variable was observed in \cite{Wi}:
\begin{equation}
\label{Will1}
\Z_{\frac{1}{n}}^{-1} \;\elaw\; n^n\,\times\, \G_{\frac{1}{n}}\,\times\,\G_{\frac{2}{n}}\,\times\,\cdots\,\times\,\G_{\frac{n-1}{n}}, \qquad n\ge 2.
\end{equation} 
A further finite factorization of $\Z_\a$ for $\a$ rational has been obtained in Formula (2.4) of \cite{TSEJP}, and reads as follows. 
\begin{equation}
\label{BetaFin}
\Z_{\frac{p}{n}}^{-p} \;\elaw \; \frac{n^n}{p^p (n-p)^{n-p}}\; \L^{n-p}\; \times\; \prod_{j=0}^{p-1} \lpa \prod_{i=q_j+1}^{q_{j+1}-1} \B_{\frac{i}{n},\frac{i-j}{n-p} - \frac{i}{n}}\rpa\; \times\; \prod_{j=1}^{p-1}\; \B_{\frac{q_j}{n},\frac{j}{p} - \frac{q_j}{n}}
\end{equation}
for every $n >p \ge 1$, where we have set $q_0 = 0, q_p = n$ and $q_j \, =\, \sup\{i\ge 1, \, ip < jn\}$ for all $j = 1, \ldots p-1.$ We refer to the paragraph before Theorem 1 in \cite{TSEJP} for more detail on this notation. 

For $\K_{\frac{p}{n}}$ and $\X_{\frac{p}{n}}$ we can obtain a finite factorization in terms of Beta random variables only, as a simple consequence of (\ref{BetaFin}). These factorizations are actually consequences of the more general Theorem 2.3 in \cite{MPZ} and Theorem 3.1 in \cite{MP}. We omit the proof.

\begin{proposition}
\label{Penson}
With the above notation, for every $n >p\ge 1$ one has
$$\K_{\frac{p}{n}}^{-p} \;\elaw \; \frac{n^n}{p^p (n-p)^{n-p}}\;\prod_{j=0}^{p-1} \lpa \prod_{i=q_j+1}^{q_{j+1}-1} \B_{\frac{i}{n},\frac{i-j}{n-p} - \frac{i}{n}}\rpa\; \times\; \prod_{j=1}^{p-1}\; \B_{\frac{q_j}{n},\frac{j}{p} - \frac{q_j}{n}}$$
and
$$\X_{\frac{p}{n}}^{-p} \;\elaw \; \frac{n^n}{p^p (n-p)^{n-p}}\;\prod_{j=0}^{p-1} \lpa \prod_{i=q_j+1}^{q_{j+1}-1} \B_{\frac{i}{n},\frac{i-j+1}{n-p} - \frac{i}{n}}\rpa\; \times\; \prod_{j=1}^{p-1}\; \B_{\frac{q_j}{n},\frac{j}{p} - \frac{q_j}{n}}.$$
\end{proposition}

\begin{remark}
\label{Pensonbis}
{\em 
(a) For $p=1,$ the above factorizations simplify into
$$\K_{\frac{1}{n}}^{-1} \;\elaw\; \frac{n^n}{(n-1)^{n-1}} \,\times\, \B_{\frac{1}{n},\frac{1}{n(n-1)}}\,\times\, \B_{\frac{2}{n},\frac{2}{n(n-1)}}\,\times\,\cdots\,\times\,  \B_{\frac{n-1}{n},\frac{1}{n}}$$
and
$$\X_{\frac{1}{n}}^{-1} \;\elaw\; \frac{n^n}{(n-1)^{n-1}} \,\times\, \B_{\frac{1}{n},\frac{n+1}{n(n-1)}}\,\times\, \B_{\frac{2}{n},\frac{n+2}{n(n-1)}}\,\times\,\cdots\,\times\,  \B_{\frac{n-1}{n},\frac{2n-1}{n(n-1)}}.$$
By the main result of \cite{BJTP}, they hence directly show that the law of $\K_{\frac{1}{n}}$ resp. $\X_{\frac{1}{n}}$ is a GGC. These Beta factorizations should also be compared to the free factorizations for $\K_{\frac{1}{n}}^{-1}$ and $\X_{\frac{1}{n}}^{-1}$ mentioned in Remark \ref{KantFuss} (a) and (b).

\vspace{2mm}

(b) In Lemma 2 of \cite{BS}, an infinite factorization of $\Z_\a^{-1}$ has also been derived in terms of Beta random variables with the help of Malmsten's formula for the Gamma function. Using this result, Corollary 1.5 in \cite{HK} and the factorization
$$\Ga_2^\beta\;\elaw\;\Gamma(\beta +2)\,\times\,\prod_{n\ge 0} \lpa\frac{n+2+\beta}{n+2}\rpa \B_{\frac{n+2}{\beta}, 1}$$
for every $\beta > 0,$ which is obtained similarly as Lemma 3 in \cite{BS}, one could be tempted to derive an infinite factorization of $\X_{\a}^{-1}$ in terms of Beta random variables for the values $\a\in (0,1)$ corresponding to the GGC property. If we try to do as in Proposition \ref{Penson}, this amounts to find factorizations of the type
$$\B_{\a + n\a, 1-\a}\;\elaw\;\B_{\frac{\a(n+2)}{1-\a}, 1}\,\times\,\B_{a_n,b_n}$$
for some $a_n,b_n >0.$ However, it can be shown that such a factorization is never possible. The existence of a suitable multiplicative factorization of $\X_\a$ which would characterize its GGC property is an open question.}
\end{remark}

In the following proposition we briefly mention a connection between $\K_{\frac{1}{n}}, \X_{\frac{1}{n}}$ and two random Vandermonde determinants, which is similar to the observations made in Section 2 of \cite{Y}. We use the notation
$$\cV(z_1,\ldots, z_n)\; =\; \prod_{1\le i<j\le n} (z_j -z_i)$$
for the Vandermonde determinant of $n$ complex numbers $z_1,\ldots, z_n.$ Let us also consider the random variable 
$$\Rr_n\; =\; \lpa\U_1\times \U_2^2 \times \cdots\times\U_{n}^{n}\rpa^{\frac{1}{n(n+1)}},$$ 
where $(\U_1,\ldots,\U_{n})$ is a sample of size $n$ of the uniform random variable on $(0,1).$   

\begin{proposition}
\label{Vander}
For every $n\ge 2,$ let $(\Tt_1,\ldots, \Tt_n)$ resp. $(\Ttt_1,\ldots, \Ttt_n)$ be a sample of size $n$ of the uniform random variable on the unit circle resp. the uniform random variable on the circle of independent random radius $\Rr_{n-1}.$ One has the identities
$$\lacc\begin{array}{lll}
 \vert \cV(\Tt_1,\ldots, \Tt_n)\vert^2  & \elaw & \K^{-1}_{\frac{1}{2}}\,\times\,\K^{-1}_{\frac{1}{3}}\,\times\,\cdots\,\times\,\K^{-1}_{\frac{1}{n}} \\
 \vert\cV(\Ttt_1,\ldots, \Ttt_n)\vert^2  & \elaw & \X^{-1}_{\frac{1}{2}}\,\times\,\X^{-1}_{\frac{1}{3}}\,\times\, \cdots\,\times\,\X^{-1}_{\frac{1}{n}}.
\end{array}
\right.
$$
\end{proposition}

\proof
To obtain the first identity, we appeal to the trigonometric version of Selberg's integral formula - see e.g. Remark 8.7.1 in \cite{AAR}, which yields
$$\esp[\vert\cV(\Tt_1,\ldots, \Tt_n)\vert^{2s}]\; =\; \frac{\Ga(1+ns)}{\Ga(1+s)^n}\; =\; \prod_{k=2}^n \lpa \frac{\Ga(1+ks)}{\Ga(1+(k-1)s)\Ga(1+s)}\rpa\; =\;  \prod_{k=2}^n \esp[\K_{\frac{1}{k}}^{-s}]$$ 
for every $s \ge 0,$ where the third equality follows at once from (\ref{Kant2}) and (\ref{FracSta}). The result follows then by Mellin inversion. The second identity is a consequence of the first one, the fact that $\vert\cV(rz_1,\ldots, rz_n)\vert^2 = r^{n(n-1)}\vert\cV(z_1,\ldots, z_n)\vert^2$ for every $r > 0$ and $z_1,\ldots,z_n\in\CC,$ and (\ref{Kant1}). 

\endproof

\begin{remark}
\label{Vanderbis}
{\em (a) If $(\Nn_1,\ldots, \Nn_n)$ is a sample of size $n$ of the standard Gaussian random variable, the Dyson-Mehta's integral formula - see e.g. Corollary 8.2.3 in \cite{AAR} - implies at once the identity 
$$\cV(\Nn_1,\ldots, \Nn_n)^2  \; \elaw \; \Z^{-1}_{\frac{1}{2}}\,\times\,\Z^{-1}_{\frac{1}{3}}\,\times\,\cdots\,\times\,\Z^{-1}_{\frac{1}{n}},$$
which is given in Proposition 3 of \cite{Y}. Observe in passing that the case $n=2$ amounts to the standard $\chi_2-$identity $\Nn_1^2\elaw  2\G_{\frac{1}{2}}.$ By (\ref{Will1}), $\cV(\Nn_1,\ldots, \Nn_n)^2$ is distributed as a finite independent product of Gamma random variables and is hence ID - see Example 5.6.3 in \cite{B}. Moreover, Theorem 1.3 in \cite{PB} and Theorem 5.1.1 in \cite{B} imply that $\vert\cV(\Nn_1,\ldots, \Nn_n)\vert$ is also ID for every $n = 4p$ or $n=4p+1.$ Since $\vert\cV(\Nn_1,\ldots, \Nn_n)\vert$ is clearly not ID for $n =2$ - see e.g. 4.5.IV in \cite{B}, one may wonder if this negative property does not hold true for every $n= 4p+2$ or $n=4p+3.$ The infinite divisibility of $\cV(\Nn_1,\ldots, \Nn_n)$ on the line seems also an open question. The logarithmic infinite divisibility of $\vert\cV(\Nn_1,\ldots, \Nn_n)\vert,$ which is easily established with explicit L\'evy-Khintchine exponent, is discussed in Section 3 of \cite{Y}.

\vspace{2mm}

(b) Setting  $\cV_n(a,b)$ for the Vandermonde determinant of $n$ independent copies of $\B_{a,b},$ a combination of the true Selberg's integral formula - see e.g.  Theorems 8.1.1 in \cite{AAR} - and Gauss's multiplication formula implies easily that $\cV_n(a,b)^2\elaw\cV_n(b,a)^2$ has a law of the type $G^{(N,N)}$ studied in Section 6 of \cite{Dufresne}, with $N = 3n(n-1)/2.$ More precisely, one has
\begin{equation}
\label{Sel}
\esp[\cV_n(a,b)^{2s}]\; =\; \frac{(a_1)_s\cdots (a_N)_s}{(b_1)_s\cdots (b_N)_s}
\end{equation}
for every $s\ge 0,$ with explicit parameters $a_i, b_i$ depending on $a$ and $b.$ For $n = 2$ and $b\ge a$ this yields the curious factorization 
$$\cV_2(a,b)^2\; \elaw\; \B_{b,a}\,\times\,\B_{\frac{1}{2},\frac{a+b}{2}}\,\times\,\B_{a, \frac{b-a}{2}}.$$
However, it does not seem that such simple Beta factorizations always exist for $n\ge 3.$ See (6) in \cite{Y} for a related identity, and also \cite{Os} for another point of view on (\ref{Sel}), where $s$ is interpreted as a parameter of a so-called Barnes Beta distribution.

\vspace{2mm}

(c) Another consequence of Proposition \ref{Vander} and Remark \ref{Pensonbis} (a) is the cyclic identity
$$ \vert \cV(\Tt_1,\ldots, \Tt_n)\vert^2\; \elaw\; n^n \, \B_{\frac{1}{n},\frac{n-1}{n}}\,\times\, \B_{\frac{2}{n},\frac{n-2}{n}}\,\times\,\cdots\,\times\,  \B_{\frac{n-1}{n},\frac{1}{n}}.$$
Let us finally mention the convergence in law
$$ \vert \cV(\Tt_1,\ldots, \Tt_n)\vert^{\frac{2}{n}}\;\claw\; e^\gamma\,\L$$
where $\gamma = -\Ga'(1)$ is Euler's constant, a simple consequence of Remark 8.7.1 in \cite{AAR}.}
\end{remark}

\subsection{An identity for the Beta-Gamma algebra} In this paragraph we prove a general identity in law which applies to the case $n=3$ in the factorizations of Proposition \ref{Penson}, and which can be viewed as a further instance of the so-called Beta-Gamma algebra - see \cite{Dufresne} and the references therein. We use the standard notation for the size-bias $X^{(t)}$ of real order $t$ of a positive random variable $X,$ that is
$$\esp[f(X^{(t)})]\; =\; \frac{\esp[X^t f(X)]}{\esp[X^t]}$$
for every $f$ bounded measurable, as soon as $\esp[X^t] < \infty.$

\begin{proposition}
\label{BG4}
For every $a,b,c, d >0$ with $a < c+d,$ one has
$$\frac{1}{\B_{a,b}\B_{c,d}}\, -\, 1\;\elaw\;\frac{\G_{b+d}}{\G_c}\,\times\,\lpa\frac{1}{1 - \B_{b,d}\B_{c+d-a,a+b}}\rpa^{(b+d)}.$$
\end{proposition}

\proof
 A direct computation using Euler's integral formula for the generalized hypergeometric functions - see e.g. (2.2.2) in \cite{AAR} - yields
$$\esp\lcr\lpa \frac{1}{\B_{a,b}\B_{c,d}}\, -\, 1\rpa^s\rcr\; =\; \esp[\B_{a,b}^{-s}]\;\esp[\B_{c,d}^{-s}]\;\pFq{3}{2}{-s,,a-s,,c-s}{a+b-s,,c+d-s}{1}$$
where we have supposed $-b-d < s < \min(a,c),$ so that the right-hand side is finite. We next appeal to Thomae's formula:
$$\pFq{3}{2}{a_1,a_2,a_3}{b_1,b_2}{1}\; =\;\frac{\Ga(b_1)\Ga(b_2)\Ga(c_1)}{\Ga(a_1)\Ga(c_1+a_2)\Ga(c_1+a_3)}\; \pFq{3}{2}{b_1-a_1,b_2-a_1,,c_1}{c_1+a_2,c_1+a_3}{1}$$ 
with $c_1 = b_1 +b_2 - a_1 - a_2 - a_3,$ which is (1) in Chapter 3.2 of \cite{Bail}, and which holds true whenever all involved parameters are positive. Setting $a_1 = a-s,$ we deduce that for every $s\in (-b-d, 0)$ one has
$$\esp\lcr\lpa \frac{1}{\B_{a,b}\B_{c,d}}\, -\, 1\rpa^s\rcr\; =\; \frac{\Ga(b+d+s)\Ga(c-s)\Ga(a+b)\Ga(c+d)}{\Ga(b+d)\Ga(c)\Ga(a)\Ga(b+c+d)}\;\pFq{3}{2}{b,,c+d-a,,b+d+s}{b+d,,b+c+d}{1},$$
and the formula extends by analyticity to $s\in(-b-d, a).$ Using again Euler's formula, the right-hand side transforms into
$$\esp[\G_{b+d}^{s}]\,\esp[\G_{c}^{-s}]\,\frac{\Ga(b+d)\Ga(c+d)}{\Ga(a)\Ga(b)\Ga(d)\Ga(c+d-a)}\int_0^1\!\!\!\int_0^1\!\! t^{c+d-a-1}(1-t)^{a+b-1}u^{b-1}(1-u)^{d-1}(1-ut)^{-b-d-s}dtdu$$
and we finally recognize
$$\esp\lcr\lpa \frac{1}{\B_{a,b}\B_{c,d}}\, -\, 1\rpa^s\rcr\; =\; \esp[\G_{b+d}^{s}]\,\esp[\G_{c}^{-s}]\,\,\frac{\Ga(a+b)\Ga(c+d)}{\Ga(a)\Ga(b+c+d)}\;\esp\lcr\lpa 1- \B_{b,d}\B_{c+d-a,a+b}\rpa^{-b-d-s}\rcr$$
for every $s\in(-b-d, a),$ which implies the required identity in law. 

\endproof

\begin{remark}
\label{BG4bis}
{\em (a) Under the symmetric assumption $c<a+b,$ we obtain the identity
$$\frac{1}{\B_{a,b}\B_{c,d}}\, -\, 1\;\elaw\;\frac{\G_{b+d}}{\G_a}\,\times\,\lpa\frac{1}{1 - \B_{d,b}\B_{a+b-c,c+d}}\rpa^{(b+d)}.$$
If both assumptions $a<c+d$ and $c<a+b$ hold, we deduce, identifying the factors and remembering $\G_s^{(t)}\elaw\G_{t+s},$ the identity
$$\G_{b+c+d}\,(1-\B_{b,d}\B_{c+d-a,a+b} )\;\elaw\; \G_{a+b+d}\, (1-\B_{d,b}\B_{a+b-c,c+d} )$$
which we could not locate in the literature on the Beta-Gamma algebra, and which boils down to the elementary $\G_{c+d}\elaw\G_{a+d}\times\B_{c+d, a-c}$ when $b = 0.$ Observe on the other hand that by Proposition 4.2 (b) in \cite{Dufresne}, this identity is equivalent to
$\B_{b,d}\G_{a+b-c} + \G_{c+d}\, \elaw\, \B_{d,b}\G_{c+d-a} + \G_{a+b},$ which is easily obtained in comparing the two Laplace transforms with the help of Euler's formula (2.2.7) in \cite{AAR}.

\vspace{2mm}

(b) Combining Propositions \ref{Penson} and \ref{BG4} yields the two identities
$$\X_{\frac{1}{3}} \, -\, b_{\frac{1}{3}}\;\elaw\; \frac{4}{27}\times\frac{\G_{\frac{3}{2}}}{\G_{\frac{1}{3}}}\,\times\lpa\frac{1}{1 - \B_{\frac{5}{6},\frac{2}{3}}\B_{\frac{1}{3},\frac{3}{2}}}\rpa^{(\frac{3}{2})}$$
and
$$\K_{\frac{1}{3}} \, -\, b_{\frac{1}{3}}\;\elaw\; \frac{4}{27}\times\frac{\G_{\frac{1}{2}}}{\G_{\frac{2}{3}}}\,\times\lpa\frac{1}{1 - \B_{\frac{1}{6},\frac{1}{3}}\B_{\frac{2}{3},\frac{1}{2}}}\rpa^{(\frac{1}{2})}.$$
Observe that these represent $\X_{\frac{1}{3}} \, -\, b_{\frac{1}{3}}$ resp. $\K_{\frac{1}{3}} \, -\, b_{\frac{1}{3}}$ as an explicit $\Ga_{\frac{3}{2}}-$mixture resp. $\Ga_{\frac{1}{2}}-$mixture, in accordance with Remark \ref{Asympbi} (b).

\vspace{2mm}

(c) Iterating Euler's integral formula (2.2.2) in \cite{AAR} yields the general representation
$$\esp\lcr\lpa \frac{1}{\B_{a_1,b_1}\ldots\,\B_{a_n,b_n}}\, -\, 1\rpa^s\rcr\; =\; \esp[\B_{a_1,b_1}^{-s}]\;\cdots\;\esp[\B_{a_n,b_n}^{-s}]\;\pFq{n+1}{n}{-s,,a_1-s,\ldots,a_n-s}{a_1+b_1-s,\ldots,a_n+b_n-s}{1}$$
for $-(b_1+\ldots + b_n) < s < \min\{a_1,\ldots, a_n\}.$ It would be interesting to know if there exists some hypergeometric transformation changing the right-hand side into
$$K\;\Ga(b_1+\cdots + b_n + s)\Ga(\max\{a_1,\ldots, a_n\}-s)\;\pFq{n+1}{n}{c_1,,\ldots,,c_n,,b_1+\cdots + b_n+s}{c_1+d_1,,\ldots,,c_n+d_n}{1}$$
for some parameters $c_i, d_i >0$ and an integration constant $K.$ This would imply the identity
\begin{equation}
\label{HGN}
\frac{1}{\B_{a_1,b_1}\ldots\,\B_{a_n,b_n}}\, -\, 1\;\elaw\; \frac{\G_{b_1+\cdots + b_n}}{\G_{\max\{a_1,\ldots, a_n\}}}\,\times\,\lpa\frac{1}{1 - \B_{c_1,d_1}\ldots\,\B_{c_n,d_n}}\rpa^{(b_1+\cdots +b_n)},
\end{equation}
which would generalize that of Proposition \ref{BG4}. Observe that in the framework of Proposition \ref{Penson} we always have $b_1+\cdots + b_n = 3/2$ resp. $1/2$ for the left-hand side of (\ref{HGN}) corresponding to $\X_{\frac{p}{n+1}}^p - 1$ resp. $\K_{\frac{p}{n+1}}^p - 1.$ Unfortunately, for $n\ge 3$ we are not aware of any such hypergeometric transformation.}
\end{remark}

\subsection{Stochastic orderings}

In this paragraph we come back to certain random variables appearing in the proof of Theorem 3. We establish some comparison results for the rescaled random variables $\V_\a = a_\a \X_\a$ with support in $[1,+\infty),$ in the spirit of those in \cite{TSEJP}. For two positive random variables $X,Y$ we write $X\prst Y$ if $\pb[X\ge x]\le\pb[Y\ge x]$ for every $x\ge 0,$ and
$$X\;\prost\; Y$$
if $X\prst Y$ and there is no such $c > 1$ such that $cX\prst Y.$ The relationship $\prost$ can be viewed as an optimal stochastic order.

\begin{proposition}
\label{StoOr}
For every $0 < \beta<\a <1$ one has
$$\frac{1}{\e}\,\times\,\U\times\W\;\prost\;\V_\a^{\frac{-\a}{1-\a}}\;\prost\;\V_\beta^{\frac{-\beta}{1-\beta}}\;\prost\; \U.$$
\end{proposition}

\proof The argument is analogous to that of (1.3) in \cite{TSEJP} and relies on (3.6) therein which, in our notation, yields
$$(a_\a\K_\a)^{\frac{-\a}{1-\a}}\;\prst\;(a_\beta\K_\beta)^{\frac{-\beta}{1-\beta}}$$
whence, by (\ref{Kant1}) and direct integration, 
$$\V_\a^{\frac{-\a}{1-\a}}\;\prst\;\V_\beta^{\frac{-\beta}{1-\beta}}$$ 
for every $0<\beta<\a<1.$ Moreover, it is easy to see by Haagerup-M\"oller's evaluation of $\esp[\X_\a^s]$ and Stirling's formula that
$$\esp[\V_\beta^{\frac{-s\beta}{1-\beta}}]\;\to\;\frac{1}{1+s}\quad\mbox{as $\beta\to 0$}\qquad\mbox{and}\qquad\esp[\V_\a^{\frac{-s\a}{1-\a}}]\;\to\;\frac{s^s}{e^s\Ga(2+s)}\quad\mbox{as $\a\to 1$}$$  
for every $s > 0.$ By Proposition \ref{MellW}, we obtain
$$\frac{1}{\e}\,\times\,\U\times\W\;\prst\;\V_\a^{\frac{-\a}{1-\a}}\;\prst\;\V_\beta^{\frac{-\beta}{1-\beta}}\;\prst\; \U.$$
To conclude the proof, by the definition of $\prost$ it is enough to observe that $\pb[\W\le\e] = 1,$ a consequence of Remark \ref{prop:SFbis} (d).

\endproof

\begin{remark}
\label{StoOrbis}
{\em (a) Multiplying all factors by an independent $\G_2$ random variable and using the second identity in Proposition \ref{prop:SF} and (\ref{SF}), we immediately retrieve Theorem A in \cite{TSEJP}.

\vspace{2mm}

(b) Proposition \ref{prop:SF} implies the limits in law
$$\lpa\frac{\a}{\a -1}\rpa\log\V_\a\;\claw\; \log\U\quad\mbox{as $\a\to 0$}\qquad\mbox{and}\qquad \lpa\frac{\a}{\a -1}\rpa\log\V_\a\;\claw\; \F \, -\, 1\quad\mbox{as $\a\to 1,$}$$ 
to be compared with that of (\ref{conv T}). This shows that distributions of the free Gumbel random variable $-\log\U$ and that of the drifted exceptional 1-free stable random variable $\F -1$ can be viewed as ``log free stable'' distributions.
 
\vspace{2mm}

(c) Specifying Proposition \ref{StoOr} to $\a =1/2$ yields
$$\frac{1}{\e}\,\times\,\U\times\W\;\prost\;\B_{\frac{1}{2},\frac{3}{2}}\;\prost\; \U,$$
whose second ordering can be observed via a single intersection property of the densities - see e.g. Lemma 1.9 (a) in \cite{DJ}. We believe the above stochastic orderings between non-explicit densities are a consequence of such a single intersection property.}
\end{remark}

Our next result deals with the classical convex ordering. For two real random variables $X,Y$, we say that $Y$ dominates $X$ for the convex order and write 
$$X\;\prcvx\; Y$$ 
if $\esp[\varphi(X)] \le \esp[\varphi(Y)]$ for every convex function such that the expectations exist.

\begin{proposition}
\label{Stocx}
For every $0 < \beta<\a <1,$ one has
$$\U\;\prcvx\;(1-\beta) {\bf X}_\beta^{\frac{-\beta}{1-\beta}}\;\prcvx\; (1-\a){\bf X}_\a^{\frac{-\a}{1-\a}} \;\prcvx\;\U\,\times\,\W.$$
\end{proposition}

We omit the proof, which is analogous to that of (1.4) in \cite{TSEJP} and a consequence of (3.7) therein. By Kellerer's theorem, this result implies that for every $t\in(0,1),$ the law of $(1-t) {\bf X}_t^{\frac{-t}{1-t}}$ is the marginal distribution at time $t$ of a martingale $\{M_t, \, t\in[0,1]\}$ starting at $\U$ and ending at $\U\times\W.$ It would be interesting to have a constructive explanation of this curious martingale connecting free extreme and free stable distributions.

\subsection{The power semicircle distribution and van Dantzig's problem}

In this paragraph we consider the power semicircle distribution with density
$$h_\a(x)\; =\; \frac{\Ga(\a+1)}{\sqrt{\pi}\Ga(\a +1/2)}\, (1-x^2)^{\a -1/2}\Un_{(-1,1)}(x),$$
where $\a > -1/2$ is the index parameter. Up to affine transformation, this law can be viewed as an extension of the arcsine, uniform and semicircle distributions which correspond to $\a = 0, \a = 1/2$ and $\a =1$ respectively. It was recently studied in \cite{AP} as a non ID factor of the standard Gaussian distribution, see also the references therein for other aspects of this distribution.

The characteristic function is computed in Formula (4.7.5) of \cite{AAR} in terms of the Bessel function of the first kind $J_\a$: one has
$${\hat h_\a}(t) \; =\; \frac{\Ga(\a+1)}{(t/2)^\a}\, J_\a(t), \qquad t > 0.$$
By the Hadamard factorization - see (4.14.4) in \cite{AAR}, we obtain
$${\hat h_\a}(z) \; =\; \prod_{n\ge 1} \lpa 1- \frac{z^2}{j^2_{\a,n}}\rpa, \qquad z\in\CC,$$ 
where $0 < j_{\a,1} < j_{\a,2} < \ldots$ are the positive zeroes of $J_\a$ and the product is absolutely convergent on every compact set of $\CC.$ 

Let now $\{X_n, \, n\ge 1\}$ be an infinite sample of the Laplace distribution with density $e^{-\vert x\vert}/2$ on $\rl$ and characteristic function
$$\esp[e^{\i t X_1}]\; =\; \frac{1}{1+t^2}\cdot$$
By (4.14.3) in \cite{AAR} and Kolmogorov's one-series theorem, the random series
$$\Sigma_\a\; =\; \sum_{n\ge 1} \frac{X_n}{j_{\a,n}}$$
is a.s. convergent. Its characteristic function is
$$\esp[e^{\i t \Sigma_\a}]\; =\;\prod_{n\ge 1} \lpa 1+ \frac{t^2}{j^2_{\a,n}}\rpa^{-1}=\; \frac{1}{{\hat h_\a}(\i t)}\cdot$$
With the terminology of \cite{L}, this means that the pair 
$$\lpa {\hat h_\a}(t), \frac{1}{{\hat h_\a}(\i t)} \rpa$$
of characteristic functions is a van Dantzig pair. The case $\a = 1/2$ corresponds to the well-known pair 
$$\lpa\frac{\sin t}{t}, \frac{t}{\sinh t}\rpa$$
which is one of the starting examples of \cite{L} and, from the point of view of the Hadamard factorization, amounts to Euler's product formula for the sine - recall from (4.6.3) in \cite{AAR} that $j_{1/2,n} = n\pi.$ The case $\a =0$ is also explicitly mentioned in \cite{L} as an example pertaining to Theorem 5 therein - observe that this theorem covers actually the whole range $\a\in(-1/2,1/2).$ In general, one has ${\hat h_\a}\in \DD_1$ for all $\a > -1/2$ with the notation of \cite{L}, and our pairs can hence be viewed as further explicit examples of van Dantzig pairs corresponding to $\DD_1.$ The case $\a =1$ is particularly worth mentioning because it shows that the semicircle characteristic function belongs to a van Dantzig pair, as does the Gaussian characteristic function.

\begin{remark}
\label{VDZ} {\em (a) The random variable $\Sigma_\a$ is ID as a convolution of Laplace distributions, and is not Gaussian. Hence, by the corollary p.117 in \cite{L}, we retrieve the fact that $X_{2,1/2}$ is not ID. Unfortunately, this method does not seem to give any insight on the non ID character of $X_{\a, \rho}$ for $\a\in (1,2)$ and $\rho\neq 1/2.$\\

(b) Following the notation of \cite{L}, the characteristic function
$${\hat g_\a}(t)\; = \; \frac{{\hat h_\a}(t)}{{\hat h_\a (\i t)}}\; =\; \frac{J_\a (t)}{I_\a(t)},$$
where $I_\a$ is the modified Bessel function of the first kind, is self-reciprocal. In other words, one has
$${\hat g_\a}(t){\hat g_\a}(\i t)\; =\; 1.$$
Observe that again, the distribution corresponding to ${\hat g_\a}(t)$ is not ID.}

\end{remark}

\subsection{Further properties of whale-shaped functions}\label{sec:WS}

In this paragraph we prove five analytical properties of WS functions and densities. Those five easy pieces apply all to the densities $f_{\a},$ and have an independent interest. We restrict the study to the class ${\rm WS}_+$, the corresponding properties for ${\rm WS}_-$ being deduced at once.

\begin{proposition}
\label{Skew}
Let $f$ be a ${\rm WS}_+$ density with unique mode $M.$ Then $f$ is  perfectly skew to the right, that is 
$$f(M+x)\; >\; f(M-x)\qquad\mbox{for every $x > 0.$}$$  
\end{proposition}

\proof Let $x_0$ be the left-extremity of ${\rm Supp}\, f$ and $M = x_1 < x_2 < x_3$ be the vanishing places of the three first derivatives of $f.$ Suppose first $M - x_0 > x_2 - M.$ Taylor's formula with integral remainder implies
$$f(M+x)\, -\, f(M-x)\; =\; \int_0^x (x-t)\,\lpa f''(M+t) - f''(M-t)\rpa\, dt.$$
On the one-hand, we have $f''(M+t) - f''(M-t) >0$ for all $t > x_2 -M$ since $f''(M-t)\le 0$ for all $t\ge 0$ and $f''(M+t) > 0$ for all $t > x_2 -M.$ On the other hand, writing
$$f''(M+t)\, -\, f''(M-t)\; =\; \int_0^t (f^{(3)}(M+s) + f^{(3)}(M-s))\, ds,$$
which is valid for all $t < M -x_0,$ we also have $f''(M+t) - f''(M-t) >0$ for all $t \le x_2 -M$ since $f^{(3)}(u) > 0$ for all $u < x_3,$ by the ${\rm WS}_+$ property. Putting everything together shows $f(M+x) > f(M-x)$ for all $x >0.$ Supposing next $M-x_0\le x_2 -M,$ the proof is analogous and easier; we just need to delete the corresponding arguments for $t > x_2 - M$.  
 
\endproof

\begin{remark}
\label{Skewbis}
{\em If we denote by $M_{\a,\rho}$ the unique mode of $f_{\a,\rho},$ the function 
$$x\;\mapsto\; f_{\a,\rho}(M_{\a,\rho} +x)\, -\, f_{\a,\rho}(M_{\a,\rho} -x)$$ 
has constant and possibly zero sign on $(0,\infty)$ for $\rho = 0, 1/2, 1$ and for $\a = 1,$ as seen from the above proposition, the explicit drifted Cauchy case and the symmetric case. One might wonder if this property of perfect skewness remains true in general. The perfect skewness of classical stable densities is a challenging open problem, which had been stated in the introduction to \cite{Ha}.}
\end{remark}

\begin{proposition}
\label{3M}
Let $f$ be a ${\rm WS}_+$ density and $M,m,\mu$ be its respective mode, median and mean. Then $f$ satisfies the strict mean-median-mode inequality
$$M\; <\; m \; < \; \mu.$$  
\end{proposition}

\proof We use the same notation of the proof of the previous proposition. First, the latter clearly implies $M < m.$ To obtain the two strict inequalities together, let us now consider the function
$$g(x)\; =\; f(m+x)\; -\; f(m-x)$$
on $[0,m-x_0].$ If $m < x_2,$ then the ${\rm WS}_+$ property implies $f^{(2)}(m+x) > f^{(2)}(m) > f^{(2)}(m-x)$ for every $x\in (0,m-x_0],$ so that $g$ is strictly convex on $(0, m-x_0].$ Since $g(0) = 0, g'(0) < 0$ and $g(m-x_0) > 0,$ this shows that $g$ vanishes only once on $(0, m-x_0]$ and from below, and hence also on the whole $(0,\infty).$ If $m \ge x_2,$ then $g$ is negative on $(0, m-M]$ and strictly convex on $[m-M, m-x_0]$ and we arrive at the same conclusion. We are hence in position to apply Lemma 1.9 (a) and (a strict, easily proved version of) Theorem 1.14 in \cite{DJ}, which implies the strict mean-median-mode inequality for $f.$  

\endproof

\begin{remark}
\label{3Mbis}
{\em (a) It is well-known and can be seen e.g. from Theorem 1.7 in \cite{HK} that $\X_{\a}$ has infinite mean. Hence, in this framework the above result only reads $M<m<\infty,$ and it is readily obtained from the previous proposition. This mode-median inequality is also conjectured to hold true for classical positive stable densities. See Proposition 5 and Remark 11 (b) in \cite{TSEJP} for partial results. 

\vspace{2mm}

(b) In the relevant case $\a \in (1,2)$ it is natural to conjecture that the strict mean-median-mode inequality holds, in one or the other direction, for both free and classical stable densities. Observe that the three parameters clearly coincide for $\rho =1/2,$ whereas for $\rho = 1/\a,$ easy computations show that the mean is zero and the mode and median are positive, so that it is enough to prove $m< M.$ In general, this problem is believed to be challenging and beyond the scope of the present paper. We refer to \cite{BDG} for a series of results on this interesting question, which however do not apply to non-explicit densities.}
\end{remark}

\begin{proposition}
\label{G2}
Let $f$ be a {\em WS} density on $(0,+\infty)$ and $X$ be the corresponding random variable. Then $X$ is a $\Gamma_2-$mixture. In particular, it is {\em ID}.  
\end{proposition}

\proof
 As in Theorem 1, we need to show that $g(x) = x^{-1}f(x)$ is a CM function, in other words that $(-1)^n g^{(n)} (x) > 0$ on $(0,\infty).$  By Leibniz's formula, we first compute
$$g^{(n)}(x)\; =\; n!\,\sum_{p=0}^n \frac{(-1)^p f^{(n-p)}(x)}{(n-p)!\, x^{p+1}}\cdot$$
This implies, after some simple rearrangements,
\begin{equation}
\label{G22}
h_n'(x)\; =\; (-1)^n x^n f^{(n+1)}(x),
\end{equation}
where $h_n(x) = (-1)^n x^{n+1} g^{(n)}(x)$ has the same sign as $(-1)^n g^{(n)}(x).$ By the WS property, we see that
$$h_n(x)\; =\; n!\,\sum_{p=0}^n \frac{(-1)^{n-p} f^{(n-p)}(x)}{(n-p)!}\,x^{n-p}$$ 
is positive on $[x_n,\infty)$ since $(-1)^i f^{(i)}(x) > 0$ when $x\in(x_i, \infty)$ for all $i\ge 0.$ Moreover, it follows from (\ref{G22}) and the whale-shape that $h_n'(x) > 0$ for $x\in(0, x_{n+1}].$  It is hence enough to show that $h_n(0+) = 0$ in order to conclude the proof, because $(0,\infty) = (0, x_{n+1}]\cup [x_n,\infty).$ But the whale-shape shows again that
$$0\; \le\; (-1)^{i-1} x^i f^{(i)}(x)\; \le\; 2\, (-1)^{i-1}x^{i-1}(f^{(i-1)}(x) - f^{(i-1)}(x/2))$$
for all $x\in (0, x_1]$ and an induction on $i,$ starting from $f(0+) = 0,$ implies $(x^i f^{(i)})(0+) = 0$ for all $i\ge 0,$ so that $h_n(0+) = 0$ as well.

\endproof

\begin{remark}
\label{G2bis}
{\em (a) The WS property is not satisfied by all densities of $\Ga_2-$mixtures vanishing at zero. A simulation shows for example that the derivative of the density
$$f(x)\; =\; x(t a^2e^{-ax}\, +\, (1-t) e^{-x})$$
vanishes three times for $a=20$ and $t =4/5.$ This contrasts with the densities of $\Ga_1-$mixtures, which are characterized by their complete monotonicity - see e.g. Proposition 51.8 in \cite{Sato}.

\vspace{2mm}

(b) For a given smooth density $f$ on $(0,\infty)$ and $n\ge 0,$ let us introduce the following property: one has $f\in{\rm BS}_n$ if 
$$\lacc\begin{array}{ll}
\sharp\{x > 0, \;\; f^{(i)}(x) = 0\} \, =\, i & \mbox{for $i \le n,$}\\
\sharp\{x > 0, \;\; f^{(i)}(x) = 0\} \, =\, n & \mbox{for $i > n.$}\\
\end{array}
\right.$$
For $n\ge 1,$ this property was introduced in \cite{TSPAMS} under the less natural denomination ${\rm WBS}_{n-1}$ - see the definition therein. Clearly, one has ${\rm BS}_0 = {\rm CM}$ and ${\rm BS}_1 = {\rm WS}$ for densities on $(0,\infty).$ Since the density of $\G_t$ has $m-$th derivative
$$(-1)^m \lpa \sum_{p=0}^m \binom{m}{p}\, (1-t)_p\, x^{-p}\rpa \frac{x^{t-1} e^{-x}}{\Ga(t)}$$
on $(0,\infty),$ it is an easy exercise using Rolle's theorem and Descartes' rule of signs to show that $\G_t\in{\rm BS}_n$ for $t\in(n,n+1].$  In this respect, the class ${\rm BS}_n$ can be thought of as an extension of the densities of $\G_t$ for $t\in(n,n+1].$ Moreover, we have just seen that the set of densities of $\Ga_{n+1}-$mixtures contains the class ${\rm BS}_n$ for $n=0,1.$ We actually believe that this is true for all $n\ge 0.$ On the other hand, the class ${\rm BS}_n$ does not seem to share any interesting property related to perfect skewness, mean-median-mode inequality or infinite divisibility for $n\ge 2.$ 

\vspace{2mm}

(c) The above proposition entails that Theorem \ref{freeID} is a consequence of Theorem \ref{shapes}. On the other hand, as we saw above, the proof of Theorem \ref{freeID} also shows that $\X_{\a}$ is a $\Ga_{3/2}-$mixture for $\a\le 3/4,$ which is not a consequence of the whale-shape.}
\end{remark}

We next study the stability of the WS property under exponential tilting. Within ID densities on $\rl,$ this transformation amounts to the multiplication of the L\'evy measure by $e^{-c\vert x\vert},$ allowing one for models with finite positive moments and analogous small jumps. This is a particular instance of the general tempering transformation, where the exponential perturbation is replaced by a CM function, and we refer to \cite{Ro} for a thorough study on tempered stable densities. If we restrict to ID densities on a positive half-line, it is seen from the L\'evy-Khintchine formula that exponential tilting amounts to multiplying the density by the same $e^{-cx}$ and renormalizing. In particular, the set of densities of $\Ga_t-$mixtures with $t\in (0,2]$ is also stable under exponential tilting.
   
\begin{proposition}
\label{Tilt}
If $f\in {\rm WS}_+,$ then $e^{-x}f\in{\rm WS}_+.$ 
\end{proposition}

\proof
 
It is enough to consider the case Supp $\! f =(0,\infty).$ Set $g(x) = e^{-x}f(x).$ Considering $h_n(x) = (-1)^n e^x g^{(n)}(x)$ for each $n\ge 0,$ we have $h_{n+1} = h_n - h_n'$ and an easy induction starting from $h_0 = f$ implies
$$(-1)^{p-1} h_{n+1}^{(p)}(0+)\, >\, 0\qquad\mbox{and}\qquad h_{n+1}^{(p)}(+\infty)\, =\, 0$$
for all $n,p\ge 0.$ We will now show that $h_{n+1}^{(p)}$ vanishes once on $(0,\infty)$ for all $n,p\ge 0,$ and that the sequence $\{x_{p,n+1}, \; p\ge 0\}$ defined by $h_{n+1}^{(p)}(x_{p,n+1})=0$ is increasing. This is sufficient for our purpose, in taking $p=0.$

Consider first the case $n=0,$ with $h_1^{(p)} = f^{(p)} - f^{(p+1)}.$ It is clear that $(-1)^ph_1^{(p)}(x) > 0$ for $x\in [x_{p+1}, \infty)$ and  that $(-1)^ph_1^{(p+1)}(x) > 0$ for $x\in (0, x_{p+1}].$ Since $(-1)^ph_1^{(p)}(0+) < 0,$ this implies that $h_1^{(p)}$ vanishes once on $(0,\infty)$ for all $p\ge 0,$ and Rolle's theorem entails that the sequence $\{x_{p,1}, \; p\ge 0\}$ defined by $h_{1}^{(p)}(x_{p,1})=0$ is increasing.

The induction step is obtained analogously from $h_{n+2}^{(p)} =h_{n+1}^{(p)} - h_{n+1}^{(p+1)},$ since $(-1)^ph_{n+2}^{(p)}(0+) < 0$ and, by the induction hypothesis, $(-1)^ph_{n+2}^{(p)}(x) > 0$ for $x\in [x_{p+1,n+1}, \infty)$ and  $(-1)^ph_{n+2}^{(p+1)}(x) > 0$ for $x\in (0, x_{p+1,n+1}].$ 
\endproof

\begin{remark}
\label{Tiltbis}
{\em (a) The above proposition implies that $e^{-x}f_{\a}$, the ``tilted free positive stable density'',  is ${\rm WS}_+$ and ID. It would be interesting to know if it is also FID.

\vspace{2mm}

(b) The class ${\rm WS}_+$ is not stable under the general tempering transformation introduced in \cite{Ro}. For example, the random variable obtained from $\G_2$ in multiplying its L\'evy measure by $te^{-x}$ is easily seen to be $(1/2)\G_{2t},$ whose density belongs to ${\rm WS}_+$ only for $t\in(1/2,1].$ }
\end{remark}

\begin{proposition}
\label{limxn}
Let $f\in{\rm WS}_+$ and $\{x_n, \, n\ge 0\}$ be the vanishing places of $\{f^{(n)}, \, n\ge 0\}.$ Then $f$ is analytic on $(x_0,\infty)$ and $x_n \to \infty.$
\end{proposition}

\proof

Again we may suppose $x_0 = 0.$ If $f$ is a density, then Proposition \ref{G2} implies that $f = xg$ where $g$ is CM and hence analytic on $(0,\infty)$, so that $f$ is analytic on $(0,\infty)$ as well. If $f$ is not a density, then Proposition \ref{Tilt} shows that $g = e^{-cx} f$ is a ${\rm WS}_+$ density on $(0,\infty)$ for some normalizing $c >0,$ and $f$ inherits the analyticity of $g$ on $(0,\infty).$

The second property is an easy consequence of the first one. Let $x_\infty$ be the increasing limit of $\{x_n, \, n\ge 0\}$ and suppose $x_\infty < \infty.$ By the whale-shape, we would then have $(-1)^nf^{(n)}(x) > 0$ for $x > x_\infty,$ so that $f$ would be CM on $(x_\infty, \infty),$ and hence also on $(0,\infty)$ by Bernstein's theorem and analytic continuation, a contradiction since $f(0+) = 0.$

\endproof

\vspace{2mm}

\noindent
{\bf Acknowledgement.} The authors were all supported by a JSPS-MAEDI research program {\em Sakura}. TH was financially supported by JSPS Grant-in-Aid for Young Scientists (B) 15K17549.

\end{document}